\documentclass[11pt,a4paper]{article}
\usepackage{amsmath,amsfonts,amssymb,mathrsfs}
\usepackage{graphicx}
\usepackage{subfigure}
\usepackage{xcolor}
\usepackage{soul}
\newtheorem{theorem}{\bf Theorem}[section]
\newtheorem{lemma}[theorem]{\bf Lemma}
\newtheorem{corollary}[theorem]{\bf Corollary}
\newtheorem{definition}[theorem]{\bf Definition}
\newtheorem{remark}[theorem]{\bf Remark}
\newtheorem{proposition}[theorem]{\bf Proposition}

\def \o {{\omega}}

\def \G {{\Gamma}}

\def \K {\mathcal{K}}
\def \L {\mathscr{L}}
\def \LY {\mathscr{L}_0}
\def \F {\mathscr{F}}
\def \A {\mathscr{A}}
\def \P {{\mathscr{P}}}

\def \H0 {{H^0_r}}
\def \HH+ {{H^+_{r, \delta}}}
\def \S+ {{S^0_{r,\delta}}}

\def \R {{\mathbb {R}}}

\def \x {{\xi}}

\def \eps {{\varepsilon}}
\def \epsilon {{\varepsilon}}

\def \t {{\tau}}
\def \t {{\tau}}

\def \th {{\theta}}

\def \g {{\gamma}}

\def \phi {{\varphi}}

\def \loc {{\text{\rm loc}}}
\def \tilde {\widetilde}

\def\p{\partial}

\def \rnn {{\mathbb {R}}^{N+1}}

\usepackage[english]{babel}
\begin{document}
\title{{Sharp Estimates for Geman-Yor Processes and applications to Arithmetic Average Asian options}
}
\author{{\sc{Gennaro Cibelli}
\thanks{Dipartimento di Scienze Fisiche, Informatiche e Matematiche, Universit\`{a} di Modena e Reggio Emilia, Via
Campi 213/b, 41125 Modena (Italy). E-mail: gennaro.cibelli@unimore.it}
\sc{Sergio Polidoro}
\thanks{Dipartimento di Scienze Fisiche, Informatiche e Matematiche, Universit\`{a} di Modena e Reggio Emilia, Via
Campi 213/b, 41125 Modena (Italy). E-mail: sergio.polidoro@unimore.it}
\sc{and Francesco Rossi}
\thanks{Dipartimento di Matematica ``Tullio Levi-Civita", Universit\`{a} degli Studi di Padova, Via Trieste 63, 35121 Padova (Italy). E-mail: francesco.rossi@math.unipd.it}
}}
\date{ }
\maketitle

\footnotesize
\begin{abstract}
We prove the existence and pointwise lower and upper bounds for the fundamental solution of the degenerate second order
partial differential equation related to Geman-Yor stochastic processes, that arise in models for option pricing theory
in finance.

Lower bounds are obtained by using repeatedly an invariant Harnack inequality and by solving an associated optimal
control problem with quadratic cost. Upper bounds are obtained by  the fact that the optimal cost satisfies a specific
Hamilton-Jacobi-Bellman equation.
\end{abstract}

\normalsize\

\section{Introduction}\label{intro}
\setcounter{section}{1} \setcounter{equation}{0} \setcounter{theorem}{0}

A keystone result in the theory of parabolic partial differential equations reads as follows: if $\Gamma =
\Gamma(x,t,\xi, \tau)$ denotes the fundamental solution of an uniformly parabolic PDE
\begin{equation*} 
  \partial_{t} u (x,t) = \sum_{i,j=1}^N \partial_{x_i} \left( a_{ij}(x,t) \partial_{x_j} u(x,t) \right), \qquad (x,t)
\in \R^N \times ]0,T],
\end{equation*}
then there exist positive constants $c^-, C^-, c^+,C^+$ such that
\begin{equation} \label{e-keystone}
  \frac{c^-}{(t-\tau)^{N/2}} \exp \left(- C^- \frac{|x-\xi|^2}{t-\tau} \right) \le \Gamma(x,t,\xi,\tau) \le
  \frac{C^+}{(t-\tau)^{N/2}} \exp \left(- c^+ \frac{|x-\xi|^2}{t-\tau} \right),
\end{equation}
for every $(x,t), (\xi,\tau) \in \R^N \times ]0,T]$ with $\tau<t$. This result has been proved by Aronson \cite{Aronson}
for operators with bounded measurable coefficients $a_{ij}$, following the fundamental works of Nash \cite{Nash} and
Moser \cite{Moser3, Moser3bis}. We also refer to the article of Fabes and Strook \cite{FaSt} for divergence form
parabolic operators, and to Krylov and Safonov \cite{KrylovSafonov2} for non-divergence form operators.

The bounds \eqref{e-keystone} have been extended by many authors to subelliptic operators. We recall in
particular, the Gaussian upper bound proved by Davies in \cite{Davies}, and the upper and lower bounds due to Jerison
and S{\'a}nchez-Calle \cite{Jerison}, and to Varopoulos, Saloff-Coste and Coulhon \cite{VSC}. We also recall that
Kusuoka and Stroock in \cite{KusuokaStroock} extend \eqref{e-keystone} by probabilistic methods. In this setting, the
quantity $|x-\xi|$ appearing in \eqref{e-keystone} is replaced by the the \emph{Carnot-Carath\'eodory} distance
$d_{CC}(x, \xi)$, that is its natural counterpart in the subelliptic setting. See also \cite{barilari}. Analogous results have been proved in
\cite{DiFrancescoPolidoro, BoscainPolidoro,CintiPolidoro, CMP, LanconelliPascucciPolidoro}, where subelliptic
parabolic operators \emph{with drift} are considered. In this case, not even the {Carnot-Carath\'eodory} distance
is appropriate to express a bound of the fundamental solution. Actually, the \emph{value function} $\Psi
= \Psi(x,t, \xi, \tau)$ of a suitable optimal control problem substitutes the whole term $\frac{|x-\xi|^2}{t-\tau}$.

In this note we extend the method used in \cite{DiFrancescoPolidoro, BoscainPolidoro, CintiPolidoro, CMP,
LanconelliPascucciPolidoro} to the study of the  degenerate parabolic operator
\begin{equation}\label{e-main}
\L u := x \partial_{x} \big( a(x,y,t) x \partial_{x} u \big) + x \, b (x,y,t) \partial_{x} u + x \partial_{y}
u - \partial_{t} u,
\end{equation}
with $(x,y,t) \in \R^+ \times \R \times ]0,T]$. The interest in the operator \eqref{e-main} arises in Finance as we
consider the problem of pricing \emph{Arithmetic Average Asian Options} in the Black \& Scholes setting. We refer to
the Black \& Scholes \cite{BlackScholes} and to Merton \cite{Merton} articles for the seminal works of this theory, and
to the books by Bj\"ork \cite{Bjork-book}, Hull \cite{Hull-book} and Pascucci \cite{Pascucci-book} for its complete
treatment. Section \ref{OptionFinance} of this article describes the application of our results to the Pricing Theory
for Financial Derivatives in the Black \& Scholes setting.

The main achievements of this article are bounds analogous to \eqref{e-keystone} for the operator $\L$. Specifically, we
prove the following inequalities for the fundamental solution $\Gamma$ of $\L$
\begin{equation} \label{e-twosidedbounds}
\begin{split}
  \frac{ c_{\eps}^-}{ t^2}\exp & \left(- C^- \Psi(x ,y + \eps t, t - \eps t \right) \le
   \Gamma(x,y,t, 1,0,0) \le
   \frac{ C_{\eps}^+}{ t^2} \exp \left(- c^+ \Psi(x,y-\eps, t + \eps) \right),
\end{split}
\end{equation}
for every $(x,y,t), \in \R^+ \times \R \times ]0,T]$ with $y + \eps t < 0$, where $\eps \in (0,1)$ is arbitrary. Here
$\Psi$ is the \emph{value function} of the  following optimal control problem
\begin{eqnarray}\label{contr-ott-0}
\Psi(x, y, t) := \inf_{\omega \in L^1([0,t])} \int_0^{t}\omega^2(\tau) d\tau \quad \mbox{subject
to constraint}
\\ \nonumber
\\ \nonumber \quad \left\{
             \begin{array}{ll}
               \dot{q_1}(s)=\omega(s)q_1(s), & \hbox{$q_1(0)= x, \quad q_1 (t)=1$,} \\
               \dot{q_2}(s)=q_1(s), & \hbox{$q_2(0)= y, \quad q_2(t)=0$.}
             \end{array}
           \right.
\end{eqnarray}
In Theorem \ref{th-main} we will give the precise statement of the bounds for $\Gamma(x,y,t, \xi, \eta, \tau)$ at any
point $(x,y,t)$  belonging to a specific subset of $\R^+ \times \R  \times [0,T]$.

\medskip

To emphasize the application of our main result to the existing literature for the operator $\L$, and to the
corresponding stochastic theory, we note that \eqref{e-keystone} can be alternatively written as
\begin{equation} \label{e-keystone-Gamma}
  {k^-}\Gamma^- (x,t,\xi,\tau) \le \Gamma(x,t,\xi,\tau) \le {k^+}\Gamma^+ (x,t,\xi,\tau),
\end{equation}
where $\Gamma^\pm$ is the fundamental solution of the heat equation $\partial_t u = \mu^{\pm} \Delta u$ with singularity
at $\xi,\tau$, and the constants $k^\pm, \mu^{\pm}$ only depend on $c^\pm, C^{\pm}$. From this point of view, it would
be natural to write \eqref{e-twosidedbounds} in terms of the fundamental solution of a suitable \emph{constant
coefficients operator} analogous to $\L$. Actually, the simplest form of $\L$ appears by choosing $a \equiv 1$, and $b
\equiv 0$:
\begin{equation}\label{e-asian}
\LY u = x^2 \partial_{xx} u + x\partial_x u+ x \partial_{y} u - \partial_{t} u, \quad (x,y,t) \in \R^+ \times \R \times
]0,T].
\end{equation}
The fundamental solution $\Gamma_0$ of $\LY$ has been first written by Yor \cite{Yor5} as the \emph{transition
density} of the process $\Big(W_t, A_t\Big)_{t\ge 0}$, where $\left(W_t\right)_{t\ge 0}$ is a Wiener process and
\begin{equation}\label{eq-Yor-process}
A_t= \int_0^t \exp \big( 2 W_s \big) \, ds.
\end{equation}
As we will see in Section \ref{Mainresults} (formula \eqref{e-Yor-density}), the expression of the fundamental solution
$\Gamma_0$ of $\LY$ is quite involved, and an estimate of the form \eqref{e-keystone-Gamma} would be hard to
handle. On the other hand, our bound \eqref{e-twosidedbounds} applies in particular to $\Gamma_0$ and provides us with
explicit information about it. Moreover, several authors point out that the explicit representation of the Asian option 
prices given by Geman and Yor \cite{Yor3} is hardly numerically treatable, in particular when pricing Asian options with 
short maturities or small volatilities (see \cite{FoschiPagliaraniPascucci, Shaw1, FuMadanWang, Dufresne}). 

Concerning the operator $\L$ in its general form, we recall that existence and regularity result for the local 
transition density were established in the recent article of Lanconelli, Pagliarani and Pascucci 
\cite{LanconelliPagliaraniPascucci}, under the assumption that the coefficients $a$ and $b$ belong to some space of
H\"older continuous functions.

\medskip

A further consequence of \eqref{e-twosidedbounds} is the following result, again in the spirit of
\eqref{e-keystone-Gamma}. By applying \eqref{e-twosidedbounds} to $\Gamma$ and to the fundamental solutions
$\Gamma^\pm$ of the operators
\begin{equation} \label{e-asian-pm}
  \L^{\pm} u =  \mu^{\pm} x^2 \partial_{xx} u + x\p_x u+ x \partial_{y} u - \partial_{t} u,
  \quad (x,y,t) \in \R^+ \times \R \times ]0,T],
\end{equation}
we obtain
\begin{equation*} 
\begin{split}
  {k^-}\Gamma^- \big(x,y + \eps (t+1), t - \eps(t+1)\big) \le & \ \Gamma(x,y,t,1,0,0) \\
  \le & \ {k^+}\Gamma^+ \left(x,y - \frac{\eps}{1-\eps}(t+1), t+\frac{\eps }{1-\eps}(t+1)\right),
\end{split}
\end{equation*}
for every $(x,y,t), \in \R^+ \times \R \times ]0,T[$ with $y + \eps (t+1) < 0$ and $t>\eps /(1-\eps)$. This is an
important theoretical result, as
it allows us to extend to $\L$ any quantitative information we know on the fundamental solution of $\L^\pm$. Clearly,
the same result holds for the densities of the respective stochastic processes. See more details in Proposition
\ref{cor-main}.

\medskip

This article is organized as follows. In Section \ref{Mainresults} we give the precise statements of our main results.
In Section \ref{OptionFinance} we explain the role which $\L$ plays in Mathematical Finance and  we give a comparison
between our bounds and similar PDE's estimates. In Section \ref{sec3} we recall known results about  the operator $\L$
defined in \eqref{e-main} and we prove a sharp Harnack inequality for it. In Section \ref{sec2}, we recall some basic
facts of stochastic processes theory, of Malliavin Calculus and we prove the existence of the density $p$ of the
stochastic process $\big(X_t, Y_t\big)_{t \ge 0}$ associated to $\L$ in \eqref{e-main}.
In Section \ref{sec4} we recall some basic tools of control theory, we solve the optimal control problem
\eqref{contr-ott-0}, and we prove the lower estimate in \eqref{e-twosidedbounds}. In Section \ref{sec6} we prove the upper
estimate in \eqref{e-twosidedbounds} and the main Theorem \ref{th-main}.

\subsection{Invariance properties and main results}\label{Mainresults}
This section contains the precise statement of our assumptions and our main results. In order to introduce the
geometrical setting useful for the study
of $\L$, we recall some properties of
$\L_0$. Monti and Pascucci observed in \cite{PascucciMonti2009} that $\LY$ is invariant with respect to the following
group
operation on $\R^+ \times \R^2$:
\begin{equation}\label{sx}
  (x_0,y_0,t_0)\circ(x,y,t)=(x_0 x,y_0+x_0 y,t_0+t).
\end{equation}
Indeed, if we set
\begin{equation}\label{eq-u-vs-v}
  v(x,y,t) = u (x_0x,y_0+x_0y,t_0+t),
\end{equation}
then $\LY v = 0$ if, and only if $\LY u = 0$. We also note that
\begin{equation} \label{e-Liegroup}
  \mathbb{G} := \left(\R^+ \times \R^2, \circ \right)
\end{equation}
is a Lie group, its identity $\mathbf{1}_\mathbb{G}$ and the inverse of $(x,y,t)$ are defined as
\begin{equation} \label{e-Lie-identity}
  \mathbf{1}_\mathbb{G} = (1,0,0), \qquad (x,y,t)^{-1} = \left(\frac{1}{x}, - \frac{y}{x}, -t \right).
\end{equation}
Then, in particular, we have
\begin{equation}\label{inv-sx}
  (x_0,y_0,t_0)^{-1} \circ(x,y,t)=\left(\frac{x}{x_0}, \frac{y-y_0}{x_0}, t-t_0 \right),
\end{equation}
so that \eqref{eq-u-vs-v} is equivalent to $u(x,y,t) = v \left(\frac{x}{x_0}, \frac{y-y_0}{x_0}, t-t_0 \right)$.

We now introduce a further notation based on the invariance properties of $\LY$ with respect to $\mathbb{G}$. As the
zero of the group $\left(\R^+ \times \R^2, \circ \right)$ is $(1,0,0)$, in the sequel we use the simplified notation
\begin{equation}\label{inv-sf}
  \G(x,y,t):=\G(x,y,t;1,0,0).
\end{equation}
Then, thanks to the invariance with respect to the left translation of $\mathbb{G}$, we have
\begin{equation*}
{x_0^2} \, \Gamma(x,y,t; x_0, y_0, t_0) = \Gamma( (x_0,y_0,t_0)^{-1} \circ(x,y,t); 1, 0,0) = \Gamma\left(\frac{x}{x_0},
\frac{y-y_0}{x_0}, t-t_0 \right).
\end{equation*}
Analogously, we denote by $\Psi(x,y,t;x_0,y_0,t_0)$ the function defined in \eqref{contr-ott-0}, with the \emph{end
point} $(1,0)$ replaced by $(x_0,y_0)$, and $t$ replaced by $t-t_0$. 
Note that, in analogy with \eqref{inv-sf}, we have
\begin{equation*}
  \Psi(x,y,t) = \Psi(x,y,t;1,0,0).
\end{equation*}

The definition of $\Psi$ is explicitly written in \eqref{prob42a} below and is well posed only when $t > t_0$ and  $y_0
> y$, otherwise problem \eqref{contr-ott-0} has no solution. In this case we agree to set $\Psi(x,y,t;x_0,y_0,t_0) =
+\infty$. The following Proposition states its invariance properties with respect to the operation on $\mathbb{G}$.

\medskip

\begin{proposition} \label{prop-Psi}
For every $(x, y, t), (x_0, y_0, t_0) \in \R^+ \times \R^2$, with $t_0 < t$ and $y_0 > y$, and for every $r >0$ we have
\begin{align} \label{tr-inv-psi}
  \Psi(x,y,t;x_0,y_0,t_0) & = \Psi\left(\tfrac {x}{x_0},\tfrac{y - y_0}{x_0},t-t_0\right); \\
\label{dil-inv-psi}
\Psi(x,y,t;x_0,y_0,t_0) & = \tfrac{1}{r}\Psi\big(x,\tfrac{y}{r},\tfrac{t}{r};
  x_0,\tfrac{y_0}{r},\tfrac{t_0}{r}\big).
\end{align}
In particular, from \eqref{tr-inv-psi} we find
\begin{equation*}
\Psi(x,y,t;x_0,y_0,t_0)=\Psi(rx,ry,t;rx_0,ry_0,t_0),
\end{equation*}
whereas, rewriting \eqref{dil-inv-psi} with $r = t-t_0$, we obtain
\begin{equation*}
\Psi(x,y,t;x_0,y_0,t_0)=\tfrac{1}{t - t_0}\Psi\left(\tfrac {x}{x_0},\tfrac{y - y_0}{(t-t_0)x_0},1\right).
\end{equation*}
\end{proposition}

\medskip

We assume the following conditions on the coefficients  of $\L$. The functions $a$ and $b$ are smooth, and there
exist two positive constants $\lambda, \Lambda$ such that
\begin{equation}\label{e-assumption1}
\begin{split}
  & | a(x,y,t) | \le \Lambda, \quad | \partial_x(x a(x,y,t)) | \le \Lambda, \quad | b(x,y,t) | \le \Lambda, \quad |
\partial_x(x b(x,y,t)) | \le \Lambda, \\
  & a(x,y,t) \ge \lambda \quad \text{for every} \quad (x,y,t) \in \R^+ \times \R^+ \times ]0,T].
\end{split}
\end{equation}

\begin{remark} \label{rem-invariance}
Unlike $\LY$, the operator $\L$ is not invariant with respect to the left translation \eqref{sx}. Indeed, as we apply
the change of variable \eqref{eq-u-vs-v} to a solution $u$ of $\L u = 0$, then $v$ is a solution of $\L_{z_0} v = 0$,
where $z_0 = (x_0, y_0,  t_0)$ and
\begin{equation}\label{Kaz0}
  \L_{z_0} v = x \partial_{x} \big( a(x_0x,y_0+x_0y,t_0+t) x \partial_{x} v \big) + x \, b (x_0x,y_0+x_0y,t_0+t)
  \partial_{x} v + x \partial_{y} v - \partial_{t} v.
\end{equation}
However, even if $\L_{z_0}$ does not agree with $\L$, it satisfies the assumption \eqref{e-assumption1} with the same
constants $\lambda$ and $\Lambda$ used for $\L$. This property will be often used in the sequel and is the basis of the
invariant nature of our bounds \eqref{e-twosidedbounds1} for the fundamental solution of $\L$.
\end{remark}

The smoothness of the coefficients $a$ and $b$ are needed to prove the existence of a fundamental solution by using
the stochastic theory (see Proposition \ref{prop-ex-fsmc} below). On the other hand, we prove upper and lower bounds for
$\Gamma$ in terms of quantities only depending on the constants $\lambda$ and $\Lambda$ appearing in
\eqref{e-assumption1}. In a future study we plan to combine the bounds \eqref{e-twosidedbounds} with, either pure PDEs
methods, or with the local results established in \cite{LanconelliPagliaraniPascucci}, to prove the existence of a
fundamental solution of $\L$ under weaker regularity assumptions on $a$ and $b$.

\medskip

We next compare our main results with the existing literature. Some results are available for the operator
$\L_0$. We  quote \cite{Yor4} for an exhaustive presentation of the topic. We mainly refer to Yor's work \cite{Yor5}
in this paper, where the author writes the density of the process \eqref{eq-Yor-process} as follows:
\begin{equation}\label{e-Yor-density}
 p(w, y, t) = \frac{e^{\frac{\pi^2}{2t}}}{\pi\sqrt{2 \pi t}}\exp\left( -\frac{1+e^{2w}}{2y}\right)
 \frac{e^w}{y^2}\psi\left(\frac{e^w}{y},t \right),
\end{equation}
where
\begin{equation}\label{psi}
  \psi\left(z,t \right)=\int_0^{\infty}e^{-\frac{\xi^2}{2 t}}e^{-z\cosh(\xi)}
  \sinh{(\xi)}\sin\left( \frac{\pi \xi}{t} \right)d \xi.
\end{equation}
Other works are due by Matsumoto, Geman and Yor \cite{Yor1, Yor2, Yor3}, Carr and Schr\"{o}der \cite{Schroder1}, Bally
and Kohatsu-Higa \cite{BallyHiga}. The fundamental solution $\Gamma_0$ of $\LY$ is
\begin{equation}\label{eq-fundsolL0}
\Gamma_0(x,y,t, x_0, y_0, t_0) =
\tfrac{1}{2xx_0}p\left(\tfrac{1}{2}\log\left(\tfrac{x_0}{x}\right),\tfrac{y_0-y}{2 x},
\tfrac{t-t_0}{2}\right).
\end{equation}

In Section \ref{sec2}, we recall some known results from the Malliavin Calculus that provide us with the existence of a
fundamental solution of $\L$ defined in \eqref{e-main}. In particular, we prove in Proposition \ref{prop-ex-fsmc} that,
if the coefficients $a$ and $b$ are smooth and satisfy suitable growth conditions, then the fundamental solution of $\L$
exists and is expressed in terms of the density of a stochastic differential equation of the form
\begin{equation} \label{eq-stoc-proc}
\begin{cases}
  & \! \! \! \! d X_t = \mu(X_t, Y_t) X_t d t + \sigma (X_t, Y_t) X_t d W_t, \\
  & \! \! \! \! d Y_t = X_t dt.
\end{cases}
\end{equation}
For this reason, in our main result we assume the existence of a fundamental solution $\Gamma$ of $\L$. We prove uniform
bounds for $\Gamma$, that only depend on the constants $\lambda$ and $\Lambda$ appearing in \eqref{e-assumption1},  and
on the $L^\infty$ norms of $a, b, \partial_x(x a)$ and $\partial_x(x b)$.

\medskip

The main result of this article is the following
\begin{theorem} \label{th-main}
Let $\Gamma$ be the fundamental solution of $\L$. Then for every $(x_0,y_0,t_0), (x,y,t) \in \R^+ \times \R \times
[0,T]$ we have
\begin{equation} \label{supp_gamma}
  \Gamma(x,y,t, x_0,y_0,t_0)=0 \qquad \forall \ (x,y,t) \in \R^+\times\R^2 \setminus
  \big\{ ]-\infty,y_0[ \times]t_0,T[ \big\}.
\end{equation}
Moreover, for arbitrary $\eps \in ]0, 1[$, there exist two positive constants $c_{\eps}^-,
C_{\eps}^+$ depending on $\eps$, on $T$ and on the operator $\L$, and two positive constants $C^-,
c^+$, only depending on the operator $\L$ such that
\begin{equation} \label{e-twosidedbounds1}
\begin{split}
  \frac{ c_{\eps}^-}{x_0^2 (t-t_0)^2} \exp & \left(- C^- \Psi(x,y+ x_0 \eps(t-t_0),t-\eps(t-t_0);x_0,y_0,t_0) \right)
\le\\
   & \quad  \quad \Gamma(x,y,t;x_0,y_0,t_0)  \le \\
   & \quad \qquad \qquad \qquad
   \frac{ C_{\eps}^+}{x_0^2 (t-t_0)^2} \exp \left(- c^+\Psi(x,y- x_0 \eps,t+ \eps;x_0,y_0,t_0) \right),
\end{split}
\end{equation}
for every $(x,y,t)\in \R^+ \times ]-\infty,y_0-x_0 \eps (t-t_0)[ \times ]t_0,T]$. Here $\Psi$ is the value function
defined in
\eqref{contr-ott-0}.
\end{theorem}

If we agree to set $\exp \left( - c^{\pm} \Psi(x,y,t;x_0,y_0,t_0) \right) = 0$ whenever
$\Psi(x,y,t;x_0,y_0,t_0)=+\infty$, then \eqref{e-twosidedbounds1} holds for every $(x_0,y_0,t_0), (x,y,t) \in \R^+
\times \R \times [0,T]$.

Clearly, the knowledge of the function $\Psi$ is crucial for the application of our Theorem \ref{th-main}.
Section \ref{sec4} of this article is devoted to the study of $\Psi$. We summarize here some of the
quantitative information about $\Psi$, that are written in terms of the function $g$ defined as follows
\begin{equation}\label{func_g}
  g(r)=\left\{
         \begin{array}{ll}
           \frac{\sinh(\sqrt{r})}{\sqrt{r}}, & \hbox{$r>0$,} \\
           1, & \hbox{$r=0$,} \\
           \frac{\sin(\sqrt{-r})}{\sqrt{-r}}, & \hbox{$-\pi^2<r<0$.}
         \end{array}
       \right.
\end{equation}

\begin{proposition} \label{prop-main}
For every $(x, y, t), (x_0, y_0, t_0) \in \R^+ \times \R^2$, with $t_0 < t$ and $y_0 > y$, 
we have
\begin{equation}\label{V(E)-prop}
\left\{
  \begin{array}{ll}
     \Psi(x_1,y_1,t_1;x_0,y_0,t_0) = E(t_1-t_0)+\frac{4(x_1+x_0)}{y_0-y_1} - 4 \sqrt{E+\frac{4x_1x_0}{(y_0-y_1)^2}},
     & \hbox{if $E\ge -\frac{\pi^2}{t_1-t_0}$;} \\
     \Psi(x_1,y_1,t_1;x_0,y_0,t_0) = E(t_1-t_0)+\frac{4(x_1+x_0)}{y_0-y_1} + 4 \sqrt{E+\frac{4x_1x_0}{(y_0-y_1)^2}},
     & \hbox{if $ -\frac{4 \pi^2}{t_1-t_0} <E <-\frac{\pi^2}{t_1-t_0}$.}
  \end{array}
\right.
\end{equation}
where
\begin{equation} \label{propE-prop}
  E = \frac{4}{(t-t_0)^2} g^{-1}\left(\frac{y_0-y}{(t - t_0) \sqrt{x x_0}}\right).
\end{equation}
Moreover,
\begin{equation} \label{asimpt1}
  \frac{\Psi(x,y,t;x_0,y_0,t_0)}{\frac{4}{(t - t_0)}\log^2\big(\tfrac{y_0 - y}{(t - t_0)
\sqrt{xx_0}}\big)+ \frac{4(x_0+x)}{y_0-y}} \to 1,
  \qquad \text{as} \qquad \frac{y_0-y}{(t - t_0)\sqrt{x_0x}} \to + \infty;
\end{equation}
\begin{equation} \label{asimpt2}
  \frac{\Psi(x,y,t;x_0,y_0,t_0)}{ \frac{4(\sqrt{x}+\sqrt{x_0})^2}{y_0-y}-\tfrac{ 4\pi^2}{(t - t_0)}} \to 1,
  \qquad \text{as} \qquad \frac{y_0-y}{(t - t_0)\sqrt{x_0x}} \to 0.
\end{equation}
\end{proposition}
The lower bound in Theorem \ref{th-main} is based on a Harnack inequality for positive solutions of $\L u=~0$. The
repeated application of the Harnack inequality, combined with a suitable optimization procedure, provides us with the
lower bound of the fundamental solution. The proof of the upper bound in Theorem \ref{th-main} for $\Gamma$
exploits the fact that the value function $\Psi$ is a solution of the relevant Hamilton-Jacobi equation.

As a corollary of Theorem \ref{th-main}, by applying \eqref{e-twosidedbounds} to $\Gamma$ and to the fundamental
solutions $\Gamma^\pm$ of the operators \eqref{e-asian-pm}, we obtain the following result. It essentially says that
the fundamental solutions of $\L$ and $\LY$ have the same behavior.

\begin{proposition} \label{cor-main}
For every $\eps \in ]0, 1[$, there exist $\Gamma^\pm$ in the form \eqref{e-asian-pm}, and positive
constants $k^{\pm}$ such that
\begin{equation*} 
\begin{split}
  {k^-}\Gamma^- (x,y + x_0 \eps (t-t_0+1) &,t- \eps (t-t_0+1);  x_0,y_0, t_0)  \le \\
  \Gamma(x,y,t, & x_0,y_0,t_0) \le
  \\ & {k^+}\Gamma^+ \left(x,y - x_0\frac{\eps}{1-\eps}(t-t_0+1),t + \frac{\eps}{1-\eps}(t-t_0+1), x_0,y_0,t_0\right),
\end{split}
\end{equation*}
for every $(x,y,t), (x_0,y_0,t_0) \in \R^+ \times \R \times ]0,T]$ with $y + x_0 \eps (t-t_0+1) <y_0$ and $t>
t_0+\eps/(1-\eps)$.
\end{proposition}

\subsection{Applications to Finance}\label{OptionFinance}
The operator $\L$ in \eqref{e-main} plays a crucial role in Mathematical Finance, since it occurs in the classical
problem of the Pricing of Arithmetic Average Asian Option. For this reason  we briefly recall in this section some
notions and details about the classic Option Pricing Theory. We start with the introduction of some simple financial
derivatives, and after we briefly recall the Black-Sholes Option Pricing Theory.
We refer to the works of Barraquand and Pudet \cite{BarraquandPudet}, and of Barucci, Polidoro and Vespri
\cite{BarucciPolidoroVespri} for a PDE approach to the pricing problem for Asian Options.

An European Put Option is a contract that gives the owner the right to sell an asset at the expiry date $T$ and at a
prescribed price $K$. A Call Option gives him, instead, the right to buy the same asset at the date $T$ and at the
price $K$. Clearly, the value of the Option at its expiry date $T$ is given by a function $\varphi(S_T)$, where $S_t$
denotes the price of the asset at time $t$. For instance, the payoff of a call option is $\varphi_C(S_T) = \max\left(0,
S_T - K\right)$, while the payoff of a put option is $\varphi_P(S_T) = \max\left(0, K - S_T \right)$. In their
celebrated article \cite{BlackScholes}, Black \& Scholes solve the problem of finding a fair price $Z = Z_t$ for this
kind of contract, at every time $t$, with $0 \le t \le T$. They assume that the price of the underlying asset at time
$t$, that is denoted by $(S_t)_{0 \le t \le T}$, is a \emph{log-normal} stochastic process,
\begin{equation}\label{e-stock}
 S_t = S_0 \exp \left( \left(\mu - \tfrac12 \sigma^2  \right) t + \sigma W_t \right), \qquad t \in  [0,T],
\end{equation}
where $(W_t)_{t\ge 0}$ denotes a standard Wiener process, $\mu$ and $\sigma$ are given constants. They construct a
self-financing portfolio, that replicates at every time $t$ the value $(Z_t)_{0 \le t \le T}$ of the Option. The
portfolio only contains an amount of the stock $(S_t)_{0 \le t \le T}$ and an amount of a riskless \emph{bond} with
constant interest rate $r$, whose price is $B_t = B_0 \exp( rt )$.
In this setting, Black \& Scholes prove that the value $ Z_t= Z(S_t,t)$ of the Option is a solution of the Black \&
Scholes equation
\begin{equation}\label{e-BS}
 \tfrac12 \sigma^2 S^2 \frac{\partial^2 Z}{\partial S^2} + r \left( S \frac{\partial Z}{\partial S} - Z \right) +
\frac{\partial Z}{\partial t} = 0, \qquad (S,t) \in  \R^+ \times ]0,T[,
\end{equation}
with final condition $Z_T = \varphi(S_T)$. We refer to Pascucci's book \cite{Pascucci-book} for an exhaustive and
detailed description of the Black \& Scholes theory and of its recent developments.

Path dependent Options are characterized by the fact that their value also depends on some average of the past price of
the stock, that is $Z_t= Z(S_t, A_t,t)$ for $0 \le t \le T$. For instance, in an \emph{Arithmetic Average Floating
Strike Option}, the strike price of an option is computed as the average of the stock price, then its payoff is
\begin{equation}\label{e-payoffAAFLS}
 \varphi_C(S_T, A_T) = \max\left(0, S_T - \frac{1}{T}\int\nolimits_0^T S_t d t \right), \quad
 \varphi_P(S_T, A_T) = \max\left(0, \frac{1}{T}\int\nolimits_0^T S_t d t - S_T \right),
\end{equation}
while in the \emph{Arithmetic Average Fixed Strike Option} the payoff is
\begin{equation}\label{e-payoffAAFXS}
 \varphi_C(S_T, A_T) = \max\left(0, \frac{1}{T}\int\nolimits_0^T S_t d t - K \right), \quad
 \varphi_P(S_T, A_T) = \max\left(0, K - \frac{1}{T}\int\nolimits_0^T S_t d t \right).
\end{equation}
When considering \emph{Geometric Average Options}, the arithmetic average $\frac{1}{T}\int_0^T S_t \, d t $ is replaced
by $\exp \left(\frac{1}{T}\int_0^T \log(S_t) \, d t \right)$.

We can summarize all the above cases by introducing the average variable $(A_t)_{0 \le t \le T}$, defined as
\begin{equation}\label{e-average}
 A_t = \int_0^t f(S_\tau) \, d\tau, \qquad t \in  ]0,T[,
\end{equation}
for some given continuous function $f$. Following the
Black \& Scholes approach, we look for the density of the process $(S_t,A_t)_{t>0}$.
We consider the stochastic differential equation of the process $(S_t,B_t,A_t)_{t>0}$,
\begin{equation}\label{e-SDE-A}
\begin{cases}
 d S_t & \! \! \! \! \! = \mu S_t  dt + \sigma S_t  d W_t, \\
 d B_t & \! \! \! \! \! = rB_t  dt,  \\
 d A_t & \! \! \! \! \! = f\left(S_t\right) dt,
\end{cases}
\end{equation}
we construct the replicating portfolio, and we apply It\^o's formula. We obtain
\begin{equation}\label{e-ABS}
 \tfrac12 \sigma^2 S^2 \frac{\partial^2 Z}{\partial S^2} + f(S) \frac{\partial Z}{\partial A}  + r \left( S
\frac{\partial Z}{\partial S} - Z \right) + \frac{\partial Z}{\partial t} = 0
\qquad (S,A,t) \in  \R^+ \times \R^+ \times ]0,T[,
\end{equation}
with \emph{final condition} $Z_T = \varphi(S_T, A_T)$.

\smallskip

We also remind that a numerical solution of the pricing problem can be obtained by a Monte Carlo method based on the
Feynman-Ka{c} formula
\begin{equation*} 
  Z(S,A,t) = \mathbb{E}^\mathbb{Q}\left[e^{-r(T-t)}\phi(S_T,A_T)|(S_t,A_t)=(S,A)\right],
\end{equation*}
where $\mathbb{Q}$ is a measure such that the process $e^{-rt}Z_t$ is a martingale under $\mathbb{Q}$.

\medskip

When considering \emph{Geometric Average Asian Option}, we have $f(S) = \log(S)$, then the simple change of variable
$v\left(e^x, y, T-t \right) := Z(S,A,t)$ transforms the PDE \eqref{e-ABS}, with its final condition, into the
following Cauchy problem
\begin{equation*} 
\begin{cases}
 & \tfrac12 \sigma^2 \left(\frac{\partial^2 v}{\partial x^2} - \frac{\partial v}{\partial x} \right)+ x
\frac{\partial v}{\partial y}  + r \left( \frac{\partial v}{\partial x} - v \right) = \frac{\partial v}{\partial t} \\
 & v(x,y,0) = \varphi \left(e^x, y \right),
\end{cases}
\end{equation*}
which, in turns, after the change of variable $u(x,y,t) := e^{rt} v\left( \frac{\sigma}{\sqrt{2}} x + \left( \tfrac12
\sigma^2 -r \right) t, y,t\right)$, can be written as follows
\begin{equation}\label{e-CP-Kolmo}
\begin{cases}
 & \frac{\partial^2 u}{\partial x^2} + x \frac{\partial u}{\partial y} = \frac{\partial u}{\partial t} \\
 & u(x,y,0) = \varphi \Big(e^{\frac{\sigma x}{\sqrt{2}}}, y \Big).
\end{cases}
\end{equation}
In PDEs theory, the solution of \eqref{e-CP-Kolmo} is given in terms of its \emph{fundamental solution} as follows
\begin{equation}\label{e-phi*Gamma}
  u(x,y,t) = \int_{\R^2} \Gamma(x,y,t, \xi, \eta, 0) \varphi \Big(e^{\frac{\sigma \xi}{\sqrt{2}}}, \eta \Big)
  d \xi d \eta.
\end{equation}
The explicit expression of the fundamental solution $\Gamma$ for the operator in \eqref{e-CP-Kolmo} is
\begin{equation} \label{e-KolmogorovFS}
\Gamma(x,y,t, \xi, \eta, \tau) =  \dfrac{\sqrt{3}}{{2 \pi}(t- \tau)^{2}}  \exp\! \left(
\!\! - \frac{|x- \xi|^2}{4(t- \tau)} - 3 \frac{|y - \eta + \tfrac{t- \tau}{2} (x +\xi)|^2}{(t- \tau)^3} \right)
\end{equation}
if $t > \tau$, while $\Gamma(x, y, t, \xi, \eta, \tau) = 0$ if $t \le \tau$ (see \cite{LanconelliPolidoro} and the
references therein).

\medskip

The function $f(S) = S$ appears in \eqref{e-ABS} as we consider \emph{Arithmetic Average Asian Option}. In this case
the function $v(x,y,t) = e^{-r t} Z(x,y,t)$ is a solution of the following PDE with final condition
\begin{equation} \label{e-CP-Kolmo-1}
\begin{cases}
 & \tfrac12 \sigma^2 x^2 \frac{\partial^2 v}{\partial x^2} + x \frac{\partial v}{\partial y}
 + r x \frac{\partial v}{\partial x} + \frac{\partial v}{\partial t} = 0, \\
 & v(x,y,T) = \varphi \left(x, y \right).
\end{cases}
\end{equation}
This problem can be further simplified by the change of variable
\begin{equation*}\label{cambiovar}
 u(x,y,t)= x^me^{m^2t}v\left(x,\frac{2y}{\sigma^2},T-\frac{2t}{\sigma^2} \right) \qquad m=\frac{r}{\sigma^2}-\frac{1}{2}
\end{equation*}
that leads to the Cauchy problem for $\L_0$
\begin{equation}\label{cauchyproblem2}
\left\{
  \begin{array}{ll}
    \L_0 u = x^2\frac{\partial^2 u}{\partial x^2}+x\frac{\partial u}{\partial x}+x\frac{\partial u}{\partial
y}-\frac{\partial u}{\partial t}=0, &
\hbox{$(x,y,t) \in \mathbb{R}^+ \times \mathbb{R}^+ \times ]0,\frac{\sigma^2}{2}T]$;} \\
    u(x,y,0)=\varphi (x,y) & \hbox{$(x,y,t) \in \mathbb{R}^+ \times \mathbb{R}^+$,}
  \end{array}
\right.
\end{equation}
whose solution writes as
\begin{equation}\label{e-phi*Gamma_0}
  u(x,y,t) = \int_{\R^+ \times \R} \Gamma_0(x,y,t, \xi, \eta, 0) \varphi (\xi, \eta ) d \xi d \eta,
\end{equation}
with $\Gamma_0$ defined  in \eqref{eq-fundsolL0}.

\medskip

The PDE approach adopted in this work allows us to consider more general problems. Among them, we can consider an
option on a basket containing $n$ assets $S_t = \big(S^1_t, \dots, S^n_t \big)$ whose dynamic is
\begin{equation} \label{eq-dyn-n}
  d S^j_t = S^j_t \mu_j(S_t, A_t, t) + S^j_t \sum_{k=1}^n \sigma_{jk} (S_t, A_t, t) d W_t^k, \qquad j=1 , \dots, n,
\end{equation}
where $\big(W^1_t, \dots, W^n_t \big)_{t \ge 0}$ is a $n$-dimensional Wiener process and $\big(A_t \big)_{t \ge 0}$ is
an average of the assets. In particular, we can choose
\begin{equation*} 
  A^j_t = \int_0^t S^j_\tau d \tau, \quad j=1, \dots, n, \qquad \text{or} \qquad
  A_t = \sum_{j=1}^n \int_0^t S^j_\tau d \tau,
\end{equation*}
including, for instance, the following
ones
\begin{equation} \label{e-moregeneral1}
\widetilde{\L}_1 u := \sum_{j,k=1}^{n} x_j \partial_{x_j} \big( a_{jk}(x,y,t) x_k \partial_{x_k} u \big) +
\sum_{j,k=1}^{n} x_j
 b_j (x,y,t) \partial_{x_j} u +\sum_{j=1}^{n}  x_j \partial_{y_j} u - \partial_{t} u,
\end{equation}
with $(x,y,t) \in (\R^+)^n \times \R^n \times ]0,T]$, and
\begin{equation} \label{e-moregeneral2}
\widetilde{\L}_2 u := \sum_{j,k=1}^{n} x_j \partial_{x_j} \big( a_{jk}(x,y,t) x_k \partial_{x_k} u \big) +
\sum_{j,k=1}^{n} x_j
 b_j (x,y,t) \partial_{x_j} u +\sum_{j=1}^{n}  x_j \partial_{y} u - \partial_{t} u,
\end{equation}
 with $(x,y,t) \in (\R^+)^n \times \R \times ]0,T]$,  respectively. In these examples, denoting by
$\sigma(x,y,t)$ the matrix $\big( \sigma(x,y,t) \big)_{j,k=1, \dots, n}$, we have
\begin{equation*} 
  \big(a_{jk}(x,y,t)\big)_{j,k=1, \dots, n} = \tfrac{1}{2}\left[ \sigma(x,y,t) \sigma(x,y,t)^* \right].
\end{equation*}
and the coefficients $b_{ij},\ i,j=1,\ldots,n$ depend on the coefficients $\mu_1, \dots, \mu_n$ and on the derivatives
of the $a_{jk}$. In this work we focus on the simplest one-dimensional case \eqref{e-main} for the sake of simplicity.

\subsubsection{Comparison with literature}\label{comparison}
We conclude this introduction with some remarks about our bounds of the fundamental solution. We first note that the
expression of $\Gamma$ in \eqref{e-KolmogorovFS} yields much information on the solution $u$. In particular, it is a
smooth function, then $u$ is smooth as well. Moreover, \eqref{e-KolmogorovFS} gives us sufficient conditions on the
function
$\varphi$ that guarantee the convergence of the integral in \eqref{e-phi*Gamma}. It is also used to prove the uniqueness
of the solution of \eqref{e-CP-Kolmo} (see \cite{Polidoro3, DiFrancescoPascucci2, DiFrancescoPolidoro}). In the same
spirit, our Theorem \ref{th-main} gives conditions on function $\varphi$ that
guarantee the convergence of the integral in \eqref{e-phi*Gamma_0}, and the uniqueness of the solution of
\eqref{cauchyproblem2} as well.
 We also compare our result with the more recent work by Delarue and Menozzi
\cite{DelarueMenozzi}, where operators in the form
\begin{equation}\label{e-DelarueMenozzi}
  \L u := \sum_{j,k=1}^{d} a_{jk}(x,t) \partial_{x_jx_k} u +\sum_{j=1}^{nd}  F_j(x,t) \partial_{x_j} u -
\partial_{t} u
\end{equation}
are considered. Here $d, n$ are positive integers, $\big(a_{jk}(x,t)\big)_{j,k=1,\dots, d}$ is a symmetric
strictly positive matrix with bounded H\"older continuous coefficients, and $F_1, \dots, F_{nd}$ satisfy suitable
assumptions. Delarue and Menozzi prove bounds for the fundamental solution of $\L$ that, in the case $d=1$ and $n=2$,
write in terms of the function $\Gamma$ in \eqref{e-KolmogorovFS}, and that of course do not apply to $\Gamma_0$ in
\eqref{eq-fundsolL0}. The reason is that, even if $\L$ in \eqref{e-moregeneral1} or \eqref{e-moregeneral2} writes in
the form \eqref{e-DelarueMenozzi}, it does not satisfy the assumption made in \cite{DelarueMenozzi}. Indeed,
following the same notations adopted in \cite{DelarueMenozzi}, our operator $\LY$ writes as above with
$$
  F_1(x,t)=0, \qquad F_2(x,t)=x, \qquad \sigma(x,t)= \left(\begin{array}{c}
                                                             x \\
                                                             0
                                                           \end{array}\right),
$$
which are respectively uniformly Lipschitz in $t$ and $\alpha$-H\"{o}lder continuous with respect to $x$, but the matrix
$$
 \tfrac12 [\sigma \sigma^*](x,t)=\left(
                           \begin{array}{cc}
                             x^2 & 0 \\
                             0 & 0 \\
                           \end{array}
                         \right),
$$
has spectrum which cannot be included in a compact interval. On the other hand, if we apply the transformation
$y =\log(x)$ we are led to consider the function
$$
  \bar{F}_1(y,t) = -\tfrac{1}{2}, \qquad \bar{F}_2(y, t)=e^{y},
  \qquad \bar{\sigma}(y, t)=\left(\begin{array}{c}
                                           1 \\
                                           0
                                           \end{array}\right),
$$
then we lose the H\"{o}lder continuity of $\bar{F}_2$ with respect to the space variable $y$.

\section{Degenerate Hypoelliptic Operators}\label{sec3}
\setcounter{section}{2} \setcounter{equation}{0} \setcounter{theorem}{0}

In this section we recall some known results about the regularity theory of linear second order operators with
non-negative characteristic form. We then introduce Harnack type inequalities and Harnack chains.

We consider a general family of differential operators, which of course contains $\L$, but also the operators defined
in
\eqref{e-moregeneral1} and \eqref{e-moregeneral2}. We set
\begin{equation} \label{e-Hormanderopeerators}
  \widetilde \L u = \sum_{i,j=1}^m X_i(a_{i,j}(z)X_j u)+\sum_{i=1}^{m} b_i(z)X_i u + Y, \qquad Y := X_0 -
\partial_t,
\end{equation}
 The  prototypes of these operator appear when we choose $a_{ij}=\delta_{ij}$ and $b_j=0$:
\begin{equation} \label{e-Hormanderopeerator}
  \widetilde \L_0 = \sum_{k=1}^m X_k^2 + Y, \qquad Y := X_0 - \partial_t,
\end{equation}
where $X_0, X_1, \dots, X_m$ are smooth vector fields defined in some open subset $\Omega$ of $\R^n \times \R$. As usual
in the PDEs theory, we identify the \emph{directional derivatives} with their \emph{vector fields}. In general, as $m <
n$, the operator $\widetilde \L_0$ is strongly degenerate. However, it may be \emph{hypoelliptic} according to the
following definition
\begin{definition}\label{hypoellipticity}
We say that $\tilde{\L}_0$ is hypoelliptic if for every distributional solution $u$  of $\widetilde \L_0  u = f$ in
$\Omega$, we have
\begin{equation} \label{e-hypoellipticity}
  u \in C^{\infty}(\Omega) \quad \text{whenever} \quad f \in C^{\infty}(\Omega).
\end{equation}
\end{definition}
The H\"{o}rmander condition \cite{Hormander} provides us with a simple sufficient condition for the hypoellipticity of
$\widetilde \L_0$. It requires the definition of \emph{commutator} of two vector fields $W$ and $Z$, acting on $u \in
C^\infty(\Omega)$ as $[W,Z]u := W Z u - Z W u$. The notation Lie$\left\{X_1, \ldots, X_m, Y \right\}(x,t)$ denotes the vector space
generate by the vector fields $\left\{X_1, \ldots, X_m, Y \right\}$ and by their commutators. The celebrated
hypoellipticity result due to H\"{o}rmander states as follows.
\begin{theorem}[H\"{o}rmander \cite{Hormander}]\label{t-Hormander}
If
\begin{equation} \label{e-Hormandercondition}
  \text{\rm Lie}\left\{X_1, \ldots, X_m, Y \right\}(x,t)= \R^n\times \R
\end{equation}
{\emph{at every $(x,t) \in \Omega$, then $\widetilde \L_0$ is hypoelliptic.}}
\end{theorem}
Concerning the operator $\LY$ in \eqref{e-asian}, we can easily check that it satisfies the H\"{o}rmander condition
\eqref{e-Hormandercondition}. Indeed, we have
\begin{equation} \label{e-vectorfields}
  X(x,y,t)=x\partial_x \sim
  \left(\begin{array}{c}
                                 x \\
                                 0 \\
                                 0 \\
                               \end{array} \right),
  \ Y(x,y,t)=  x \partial_y-\partial_t \sim
  \left(\begin{array}{c}
                                 0 \\
                                 x \\
                                -1 \\
                                \end{array} \right),
  \   [X,Y](x,y,t)= x \partial_y \sim \left(\begin{array}{c}
              0 \\
              x \\
              0 \\
              \end{array} \right).
\end{equation}
Then, the vectors $X, Y$ and $[X,Y]$ form a basis of $\R^3$ at every point $(x,y,t) \in \R^+ \times \R^2$. By
H\"ormander's Theorem \ref{t-Hormander}, $\LY$ is hypoelliptic in $\R^+ \times \R^2$ in the sense of Definition
\ref{hypoellipticity}.
In PDE's Theory the regularity of operators satisfying H\"{o}rmander condition is strongly related to a Lie group
structure on the underlying domain. We refer to the seminal works of Folland \cite{Folland}, Folland-Stein
\cite{FollandStein}, Nagel-Stein-Wainger \cite{nagelstein}.

\smallskip

For the sake of clarity, we now recall the definition of fundamental solution for a hypoelliptic operator $\L$. With
this aim we write $\L$ in its divergence form
\begin{equation}\label{divergence_form}
  \L u =  -X^*(aXu)+(b-a)Xu+ Yu, 
\end{equation}
where
$X^* u(x,y,t) := -Xu(x,y,t)-u(x,y,t).$
\begin{definition}\label{def_fundsol}
We say that a function $\Gamma: (\R^+\times\R^2) \times (\R^+\times \R^2) \to \R$ is a \emph{fundamental solution} of
$\L$ if:
\begin{enumerate}
\item for every $(x_0,y_0,t_0) \in \R^+\times\R^2$ the function $(x,y,t) \mapsto \Gamma(x,y,t;x_0,y_0,t_0)$:
\begin{description}
\item[{\it i)}] belongs to $L^1_\loc(\R^+\times\R^2) \cap C^{\infty}(\R^+\times\R^2
\setminus\left\{(x_0,y_0,t_0)\right\} )$,
\item[{\it ii)}] it is a classical solution of $\L u = 0$ in $\R^+\times\R^2 \setminus\left\{(x_0,y_0,t_0)\right\}$;
\end{description}
\item for every $\phi \in C_b(\R^2)$ the function
$$
  u(x,y,t)=\int_{\R^+\times \R}\Gamma(x,y,t;\xi,\eta,0)\phi(\xi,\eta)d \xi\, d\eta,
$$
is a classical solution of the Cauchy problem
\begin{equation} \label{cauchyproblem}
\left\{
  \begin{array}{ll}
    \L u = 0, & \hbox{$(x,y,t) \in \mathbb{R}^+ \times \mathbb{R} \times \mathbb{R}^+$;} \\
    u(x,y,0)=\varphi (x,y) & \hbox{$(x,y,t) \in \mathbb{R}^+ \times \mathbb{R}$.}
  \end{array}
\right.
\end{equation}
\item The function $\Gamma^*(x,y,t;x_0,y_0,t_0): = \Gamma(x_0,y_0,t_0; x,y,t)$ satisfies 1. and 2. with $\L$ replaced
by its formal adjoint
\begin{equation}\label{main_adj}
\L^* v :=  - X^* \big( a X v \big) + X^* \big((b-a) v \big) - Y v.
\end{equation}
\end{enumerate}
\end{definition}

The main tool in the proof of our asymptotic estimates of the fundamental solution are the Harnack inequalities and the
\emph{Harnack chains}. In this setting a Harnack chain is defined as follows:
\begin{definition} \label{def-HC}
Let $\Omega$ be an open subset of $\rnn$. We say that a finite set $\{z_0, z_1,...,z_k \} \in
\Omega$ is a  \emph{Harnack chain} connecting $z_0$ to $z_k$ if there exist positive constants $C_1,...,C_k$ such that:
$$
  u(z_j) \leq C_j u(z_{j-1}), \qquad j=1,...,k,
$$
for every positive solution $u$ of $\widetilde{\L} u = 0$.
\end{definition}

Harnack chains have been used by several authors to prove asymptotical lower bounds of the fundamental solution of
degenerate hypoelliptic operators. See for instance \cite{VSC, Polidoro1, DiFrancescoPolidoro, BoscainPolidoro, CMP,
PascucciPolidoro6}. They have been also used to prove asymptotic estimates near the boundary for the positive solution
of Kolmogorov operators, see \cite{CNP1, CNP2}. In the above articles, Harnack chains have been constructed by
selecting
points belonging to the trajectories of $\widetilde \L$-admissible paths, which are defined as follows:

\begin{definition}
An $\widetilde \L$-admissible path with starting point $z_0$ is a solution of the following Cauchy problem
\begin{equation}\label{camminiamm}
 \dot{\gamma}(s)= \sum_{k=1}^m \omega_k(s)X_k(\gamma(s))+Y(\gamma(s)), \qquad \gamma(0)=z_0
\end{equation}
where $\omega(s)=(\omega_1(s),\dots,\omega_m(s)) \in \mathbb{R}^{m}$, $s\geq0$ and each $\omega_i(s) \in
L^1[0,+\infty[$.
\end{definition}

We next focus on the operator $\L$ in \eqref{e-main}.

\subsection{Harnack inequality and Green function for $\L$}

Our construction of Harnack chains for $\L$ is based on the following Harnack inequality. Its statement requires some
notation. For any $z_0 = (x_0, y_0, t_0) \in \R^+ \times \R^2$ and $r \in ]0,1[$, we set

\begin{equation}\label{spheres}
\begin{split}
    H_r(z_0)& = \Big\{ (x,y,t) \in \R^3 : |x - x_0| < r x_0, -r^2 < t-t_0 < 0,
    | y-y_0 + x_0(t-t_0)| < r^3x_0 \Big\}\\
    S_r(z_0) & = \Big\{ (x,y,t) \in \R^3 : |x - x_0| \le r x_0, - r^2 \le
    t- t_0 \le - \frac{r^2}2, | y-y_0 + x_0(t-t_0)| \le r^3x_0 \Big\}
\end{split}
\end{equation}

Notice that the cylinders defined in \eqref{spheres} are the most natural
geometric sets which can be defined taking into account the  invariance group \eqref{sx}
of $\L_0$. Indeed, they are obtained from $H_r(1,0,0)$ and $S_r(1,0,0)$, respectively, by using the left translation
``$\circ$'' in \eqref{sx}.
\begin{proposition}\label{p1a}
Let $z_0 \in \R^+ \times \R^2$ and $r \in ]0, 1/2]$.
 If $u$ is a positive solution of $\L u = 0$ in $H_r(z_0)$,
then
\begin{equation*}
  u(z)\le M \, u(z_0)
\end{equation*}
for every $z \in S_{\theta r}(z_0)$. The two constants $\theta \in ]0,1[$ and $M >0$ only depend on the operator $\L$.
\end{proposition}

The proof of Proposition \ref{p1a} relies on the Harnack inequality proven by Golse, Imbert, Mouhot, and Vasseur in
\cite{GIMV}. We also refer to \cite{AEP} for a geometric statement of the Harnack inequality. The operators  $\K$
considered in \cite{GIMV} and \cite{AEP} act on a function $u$ as follows
\begin{equation} \label{e-Kolmo-div}
\K u := \sum_{j,k=1}^{n} \partial_{x_j} \big( \widetilde a_{jk}(x,y,t) \partial_{x_k} u \big) +
\sum_{j,k=1}^{n} \widetilde b_j (x,y,t) \partial_{x_j} u + \sum_{j=1}^{n}  x_j \partial_{y_j} u - \partial_{t} u.
\end{equation}
Here $(x,y,t) \in \R^n \times \R^n \times \R$ and the coefficients $\widetilde a_{jk}, \widetilde b_{j}$ are bounded
measurable functions for $j,k = 1, \dots, n$. Moreover $\widetilde a_{jk}= \widetilde a_{kj}$ and
\begin{equation} \label{ellipticity-K}
  \sum_{j,k=1}^{n}\widetilde a_{jk}(x,y,t)\xi_{j}\xi_{k} \ge \lambda |\xi|^{2},\qquad \text{for every} \
  \xi \in \R^{n},\ \text{and} \ (x,y,t)\in \R^{2n+1}.
\end{equation}
Note that the main structural difference between $\L$ and $\K$ is  in that the
coefficients of $\K$ are bounded and satisfy the \emph{unform ellipticity} condition \eqref{ellipticity-K}, with respect
to the variable $x \in \R^n$. As the Harnack inequality is a \emph{local} result, we will borrow the Harnack inequality
for $\K$ for the study of the positive solutions of $\L u = 0$. For the sake of simplicity, we recall the statement of
the Harnack inequality proven in \cite{GIMV} only for $n=1$ and with a notation suitable for our operator $\L$.

{\it Let $\Omega$ be an open subset of $\R^3$. Consider the following operator}
\begin{equation} \label{e-kolmo-a}
  \K v =  \partial_{x} \left(\widetilde a(x,y,t) \partial_{x} v \right) + \widetilde b(x,y,t) \partial_{x} v + x
\partial_{y} v -
\partial_{t} v, \quad (x,y,t) \in \Omega.
\end{equation}
{\it Assume that $\widetilde a$ and $\widetilde b$ are bounded  measurable functions such that
$\inf_{\R^3} \tilde{a}(x,y,t) >0$. Let $z_0 \in \Omega, \ r \in ]0,1/2]$
be such that $H_r(z_0)\subseteq\Omega$. Then there exist two positive constants $\th$ and $M$, only depending on the
operator $\K$, such that}
\begin{equation}\label{harnack-K}
  v(z) \le M \, v(z_0), \qquad \text{for every} \ z \in S_r(z_0),
\end{equation}
{\it and for every non-negative solution $v$ of $\K v = 0$ in $\Omega$.}

\medskip

\noindent {\sc Proof of Proposition \ref{p1a}.} \ Let $u$ be a positive solution of $\L u = 0$ in $H_r(z_0)$, with $r
\in ]0, 1/2]$. We first consider the  point  $z_0 =(1,0,0)$. With the aim
to apply \eqref{harnack-K} to $u$,  we write $\L$ in the
form \eqref{e-kolmo-a} by setting
\begin{equation}\label{e1ab}
    \widetilde a(x,y,t) = x^2 a(x,y,t), \qquad
     \widetilde b(x,y,t) = x \left(b(x,y,t) - a(x,y,t) \right).
\end{equation}
 In order to deal with
\emph{bounded} coefficients $\widetilde a$ and $\widetilde b$, we modify them out of the cylinder $H_r(z_0)$ as
follows. We set
\begin{equation}\label{e2ab}
    \widetilde a(x,y,t) := \varphi^2(x) a(x,y,t), \qquad
    \widetilde b(x,y,t) := \varphi(x) \left( b(x,y,t) - a(x,y,t) \right),
\end{equation}
where
\begin{equation}\label{e1a}
    \varphi(x) = \left\{
    \begin{array}{ll} 1/2 & \text{for} \ x \in ]0,1/2]  , \\ 
    x & \text{for} \ x \in ]1/2,3/2[, \\ 
    3/2  & \text{for} \ x \in [3/2, \infty[. 
    \end{array} \right.
\end{equation}
Then, it is easy to check that our assumption \eqref{e-assumption1}  on $\L$
implies the conditions on $\K$ for the validity of \eqref{harnack-K}. In particular, our claim is proven for $z_0
=(1,0,0)$ and for every $r \in ]0, 1/2]$, since in this case $\L$ agrees with $\K$ in the cylinder $H_r(z_0)$.

\medskip

An argument similar to that used above would give the proof of Proposition \ref{p1a} with a constant $M$ that may
depend on $z_0$. In order to prove our claim as stated, with $M$ independent on $z_0$, we rely on the left translation
\eqref{sx}. As we apply the change of variable \eqref{eq-u-vs-v} to a solution $u$ of $\L u = 0$ in $H_r(z_0)$, then $v$
is a solution of $\L_{z_0} v = 0$ in $H_r(1,0,0)$ where $\L_{z_0}$ is defined in \eqref{Kaz0}.
Note that, as we have noticed in Remark \ref{rem-invariance}, $\L_{z_0}$ satisfies assumptions \eqref{e-assumption1}, with the same constants used for $\L$. In particular, the Harnack inequality
\eqref{harnack-K}
holds for $v$, and implies
\begin{equation*} 
  u(x,y,t) = v \left(\tfrac{x}{x_0}, \tfrac{y-y_0}{x_0}, t-t_0 \right) \le M \, v(1,0,0) = M \, u(x_0, y_0,t_0),
\end{equation*}
for every $x,y,t \in S_r(x_0, y_0,t_0)$. This concludes the proof. \hfill $\square$

\medskip

As a direct consequence, we obtain the following
\begin{corollary}\label{c1a}
If $u$ is a positive solution of $\L u = 0$ in $ H_r(z_0)$, where $0<r\leq 1/2$, 
 then
\begin{equation*}
  u(z)\le M \, u(z_0)
\end{equation*}
for every $z$ in the set
\begin{equation}\label{paraboloid-Ka}
\begin{split}
  \P_r(z_0) = \Big\{ (x,y,t) \in \R^3 : 0 < t_0-t \le \theta^2 r^2,
  & |x - x_0| \le (t_0-t)^\frac12 x_0, \\ &
  |y-y_0 - (t_0-t) x_0| \le (t_0-t)^\frac32 x_0 \Big\}.
\end{split}
\end{equation}
\end{corollary}

\medskip
A crucial ingredient for  the proof of our lower bound of the fundamental solution to $\L$  is the analogous lower bound of a Green function $G$ for the operator $\K$ defined in
\eqref{e-kolmo-a}. The existence of a Green function for $\K$  \eqref{e-kolmo-a} has been established by Di Francesco and Polidoro in \cite{DiFrancescoPolidoro} if the
coefficients $\widetilde a_{jk}, \widetilde b_{j}, j,k = 1, \dots, n$ are bounded, and H\"older continuous
 functions, and \eqref{ellipticity-K} is satisfied. In \cite{DiFrancescoPolidoro} it is also given a
lower bound for $G$, in terms of constants depending on the H\"older continuity of the coefficients of $\K$ (see Theorem
4.3 in \cite{DiFrancescoPolidoro}). Here we give a bound of $G$ where the constants only depend on the dimension $n$, on
the constant $\lambda$ in \eqref{ellipticity-K} and on the $L^\infty$ norm of $\widetilde a_{jk}, \widetilde b_{j}, j,k
= 1, \dots, n$. We rely on the method used in \cite{DiFrancescoPolidoro} and on the upper and lower bounds proven by
Lanconelli, Pascucci and Polidoro in \cite{LanconelliPascucciPolidoro} (see also \cite{LanconelliPascucci}).

We next recall the statement Theorem 1.3 in \cite{LanconelliPascucciPolidoro} with the notation used here for the
operator $\K$. Here $\Gamma_\K$ denotes the fundamental solution to $\K$, while  $\Gamma_\K^\mu$ is the fundamental
solution to the \emph{constant coefficients operator}
\begin{equation*} 
 \K^{\mu}:= \mu \sum_{j=1}^{n} \p_{x_{j}x_{j}} + \sum_{j=1}^{d}b_{j}x_j \partial_{y_j} -  \partial_t.
\end{equation*}
\emph{Assume that the coefficients $\widetilde a_{jk}, \widetilde b_{j}, j,k = 1, \dots, n$ of the operator $\K$ are
bounded measurable functions and  that \eqref{ellipticity-K} is satisfied. Let $I = ]T_0, T_1[$ be a bounded interval.
Then,  there exist four positive constants $\mu^+, \mu^-, C^+, C^-$ such that}
\begin{equation} \label{eq-bounds}
 C^- \Gamma_\K^{\mu^-}(x,y,t, \xi, \eta,\tau) \le \Gamma_\K(x,y,t, \xi, \eta,\tau)
 \le C^+ \Gamma_\K^{\mu^+} (x,y,t, \xi, \eta,\tau),
\end{equation}
\emph{for every $(x,y,t), (\xi, \eta,\tau) \in \R^{2n+1}$ with $T_0 < \tau < t < T_1$. The constants $\mu^-,
\mu^+$ depend only on $n$ and $\L$, while $C^-, C^+$ also depend on $T_1-T_0$.}

\medskip

We recall that the explicit expression of $\Gamma_\K^{\mu^{\pm}}$ is known (see, for instance \cite{Hormander} and
\cite{LanconelliPolidoro}):
\begin{equation} \label{eq-L-lambda}
 \Gamma_\K^{\mu}(x,y,t, \xi, \eta,\tau) =  \frac{3^{n/2}}{  (2 \pi \mu)^n (t-\tau)^{2n}}
    \exp\bigg(-\frac{1}{4 \mu} \bigg(\frac{|x-\xi|^2}{t-\tau}  +
    12 \frac{\left|y - \eta + (t-\tau) (x + \xi) /2 \right|^2}{(t-\tau)^3} \bigg) \bigg),
\end{equation}
for every $\tau<t$ and $(x,y), (\xi, \eta) \in\R^{2n}$. We also recall that $\Gamma_\K^{\mu^{\pm}}$ are homogeneous
of degree $-4n$ with respect to the dilation $(x,y,t) \mapsto (r x, r^3 y, r^2 t)$, that is
\begin{equation} \label{eq-L-dilation}
 \Gamma_\K^{\mu^{\pm}}(r x,r^3 y,r^2 t, r \xi, r^3 \eta,r^2 \tau) =  \frac{1}{ r^{4n}} \,
 \Gamma_\K^{\mu^{\pm}}(x,y,t, \xi, \eta,\tau),
\end{equation}
for every $(x,y,t; \xi, \eta, \tau) \in R^6$ and for every positive $r$.

We next recall the method used in \cite{DiFrancescoPolidoro} to prove a lower bound for the Green function, in
order to remove the H\"older regularity assumption made on the coefficients of the operator. With this aim, we
introduce here a simplified notation useful for our purpose. We first define a cylinder analogous to $H_r(z_0)$,
centered at $z_0 = (1,0,0)$. For any $r, \delta \in ]0,1/2]$, we set
\begin{equation}\label{eq-cylinder}
\begin{split}
   \H0 (1,0,0) & = \Big\{ (x,y,t) \in \R^3 : \tfrac{(x - 1)^2}{r^2} + \tfrac{|x - 1|}{r}+
   \tfrac{(y + t)^2}{r^6} < 1, 0 < t  <{r^2}  \Big\}, \\
   \S+ (1,0,0) & = \Big\{ (x,y,t) \in \R^3 : \tfrac{(x - 1)^2}{r^2} + \tfrac{|x - 1|}{r}+
   \tfrac{y^2}{r^6} \le \delta, t = 0  \Big\}.
\end{split}
\end{equation}

\medskip
Note  that  $\H0 (1,0,0) \subset \big\{1 - r <  x < 1+r\big\}$. In particular, if we define $\widetilde a$ and
$\widetilde b$ according to \eqref{e2ab}  and \eqref{e1a}, then $\K$ agrees with $\L$ in the cylinder $\H0 (1,0,0)$.
Also note that the geometry of $\H0 (1,0,0)$ is more complicated than the one of $\H0 (1,0,0)$. The advantage of
this fact is that the the Dirichlet problem for $\K$ in \eqref{e-kolmo-a} is well posed in $\H0 (1,0,0)$ .

In Section 4 of \cite{DiFrancescoPolidoro} it is proven the existence of a \emph{Green function} $G_r:
\overline{\H0 (1,0,0)} \times \H0 (1,0,0) \to [0, + \infty[$ with the following property: for every $f \in C_0^\infty
(\H0 (1,0,0))$, the function
\begin{equation} \label{eq-v-Green1}
  v_r (x,y,t) := \int_{\H0 (1,0,0)} G_r(x,y,t; \xi, \eta, \tau) f (\xi, \eta, \tau) d\x \, d \eta \, d \tau,
\end{equation}
is a classical solution of the Dirichlet problem
\begin{equation} \label{eq-v-Green2}
\begin{cases}
  & \L u =  - f \quad \text{in} \ \H0 (1,0,0), \\
  & \ \ \ \ u = 0 \quad \text{in} \ \partial (\H0 (1,0,0)) \cap \big\{t < T\big\}.
\end{cases}
\end{equation}
The Green function $G_r$ for the cylinder $\H0 (1,0,0)$ is defined in \cite{DiFrancescoPolidoro} as follows:
\begin{equation}\label{eq_Green_Kolm}
G_r(x,y,t; \xi, \eta, \tau) = \Gamma_{\K}(x,y,t, \xi, \eta, \tau)-h_r(x,y,t, \xi, \eta, \tau),
\end{equation}
where  $h_r(x,y,t;\xi,\eta,\tau)$ is the solution to the Dirichlet problem:
\begin{equation} \label{eq-v-Green3}
\begin{cases}
  & \L u = 0 \quad \text{in} \ \H0 (1,0,0), \\
  & \ \ \ \ u = \Gamma_{\mathcal{K}}(x,y,t;\xi,\eta,\tau) \quad \text{in}
  \ \partial (\H0 (1,0,0)) \cap \big\{t < T\big\}.
\end{cases}
\end{equation}
The following result will be needed in the proof of the lower bound of the fundamental solution.

\begin{lemma} \label{lem-bd-Green} There exist two positive constants $\kappa$ and $\varrho$, only depending on the
$L^\infty$ norms of $\widetilde a, \widetilde b$, and on $\inf \widetilde a$, such that
  \begin{equation*} 
  G_{r}(1,- s, s; 1,0,0) \ge \frac{\kappa}{s^2}, \quad \text{for every} \quad  s \in ]0,  \varrho r^2[.
\end{equation*}
\end{lemma}

\medskip

\noindent{\sc Proof.} Choose any $r, \delta \in ]0,1/2]$, and consider the compact set
$$
    M_r(1,0,0) := 
    \overline{\partial (\H0 (1,0,0))\cap \{0<t<T\}} \times \S+ (1,0,0).
$$
Let $(\xi,\eta,\tau)$ be a point of $\S+ (1,0,0)$, and let $h_r$ be the solution to \eqref{eq-v-Green3}. By
the strong maximum principle we have that $h_r \ge 0$ and
$$
  \max_{(x,y,t) \in \overline {\H0 (1,0,0)}} h_r (x,y,t) = \max_{(x,y,t) \in \overline {\partial (\H0 (1,0,0))\cap
\{0<t<T\}}} \Gamma_{\mathcal{K}}(x,y,t;\xi,\eta,\tau).
$$
Then, by using \eqref{eq-bounds} in the above inequality, We find that
$$
  \max_{(x,y,t) \in \overline {\H0 (1,0,0)}} h_r (x,y,t) \le  \widetilde  \kappa_r, \qquad
  \widetilde  \kappa_r :=  C^+ \max_{(x,y,t) \in M_r(1,0,0)} \Gamma_{\mathcal{K}}^{\mu^+}(x,y,t;\xi,\eta,\tau).
$$
We also note that $\widetilde  \kappa_r = \frac{\widetilde  \kappa_1}{r^4}$ because of \eqref{eq-L-dilation}. As a
consequence of the above inequalities, of \eqref{eq-bounds} and of the definition \eqref{eq_Green_Kolm} of $G_r$ we then
find
\begin{equation*} 
G_r(x,y,t; \xi, \eta, \tau) \ge C^- \Gamma_\K^{\mu^-}(x,y,t, \xi, \eta,\tau) - \frac{\widetilde  \kappa_1}{r^4}
\end{equation*}
for every $(\xi,\eta,\tau) \in \S+ (1,0,0)$ and $(x,y,t) \in \H0 (1,0,0)$. In particular,
  \begin{equation} \label{eq-bd-Green}
  G_{r}(1,- s, s; 1,0,0) \ge C^- \Gamma_\K^{\mu^-}(1,- s, s; 1,0,0) - \frac{\widetilde  \kappa_1}{r^4} =
\frac{C^{-}\sqrt{3}}{2 \pi \mu^- s^2 } - \frac{\widetilde  \kappa_1}{r^4},
\end{equation}
for every $s \in ]0, r^2[$. We eventually choose a positive $\kappa$ such that $\kappa < \frac{C^{-}\sqrt{3}}{2 \pi
\mu^-}$ and we conclude that there exist a positive $\varrho$ such that
  \begin{equation*} 
  \frac{C^{-}\sqrt{3}}{4 \pi \mu^- s^2 } - \frac{\widetilde  \kappa_1}{r^4} >
\frac{\kappa}{s^2}, \quad \text{for every} \quad s \in ]0, \varrho r^2[.
\end{equation*}

This inequality and \eqref{eq-bd-Green} conclude the proof. \hfill $\Box$

\subsection{Harnack chains for $\L$}

Any $\L$-admissible path $\g(s)=(x(s),y(s),t(s))$ for $\LY$ is
the solution of the Cauchy problem
\begin{equation}\label{paramet-gamma}
  \left\{
    \begin{array}{ll}
      \dot{x}(s)=\o(s)x(s) &\hbox{$x(0)=x_0$,}\\
      \dot{y}(s)=x(s) &\hbox{$y(0)=y_0$,}\\
     \dot{t}(s)=-1, & \hbox{$ t(0)=t_0$,}
    \end{array}
  \right.
\end{equation}
where $\o \in L^1([0,t_0-t])$.
In this setting, we refer to the function $\omega$ as the \emph{control} of the problem
\eqref{paramet-gamma}. We introduce now a standard definition from control theory, see \cite{agra-book}:

\medskip

\begin{definition} (Attainable set).
For every $z_0 \in \Omega \subseteq \mathbb{R}^{3}$  the
\emph{attainable set}  $\A_{z_0}$ from $z_0$ in $\Omega$ is
\begin{align}\label{insprop}
  \A_{z_0}=\big\{z \in \Omega \mid &\mbox{there exists a time $\bar{t} \in \R^+$ and an $\L$-admissible path} \nonumber
\\
   & \qquad \qquad \qquad \qquad \ \gamma: \
[0,\bar{t}]\rightarrow \Omega \ s.t \ z_0= \g(0),\ z = \gamma(\bar{t}) \big\}.
\end{align}
\end{definition}

\medskip

\begin{proposition}\label{supp-density}
  For every $(x_0,y_0,t_0) \in \mathbb{R}^+\times\mathbb{R} \times ]T_0,T_1[$ it holds:
 \begin{equation}\label{propagazione}
   \A_{(x_0,y_0,t_0)}=
   ]0,+\infty[\times]y_0,+\infty[\times ]T_0,t_0[.
\end{equation}
\end{proposition}

\medskip

\noindent {\sc Proof.} \ From \eqref{paramet-gamma} it plainly follows that
$$
  \A_{(x_0,y_0,t_0)}(\mathbb{R}^+\times\mathbb{R} \times]T_0,T_1[) \subseteq
  [0,+\infty[\times]y_0,+\infty[\times ]T_0,t_0[.
$$
The opposite inclusion will follow from the results given in Section 4.2, where we
exhibit an $\L$-admissible path steering $(x_0,y_0,t_0)$ to any given point $(x,y,t) \in ]0,+\infty[ \times
]y_0,+\infty[ \times ]T_0,t_0[$.
\hfill $\square$

\medskip

The following result provides us with a bound of any positive solution $u$ of $\L u = 0$ at the end point
$\g(t_0-t)$ of an $\L$-admissible path $\gamma$.

\begin{proposition}\label{p2a}
There exist four positive constants $\theta, h, \beta$ and $M$, with $\theta < 1$ and $M>0$ , only depending on the
operator $\L$ such that the following property holds.

Let $T_0<t<t_0<T_1$ be fixed. Fix $(x_0,y_0)$ and let $\o\in L^1([t,t_0],\R)$ be a control, with $\g : [t,t_0] \to
\R^3$ the corresponding $\L$-admissible path of
\eqref{paramet-gamma} starting from $(x_0,y_0,t_0)$. Denote by $(x,y,t)=\g(t_0)$ its end-point. Then, for every
positive solution $u: \R^+ \times \R \times ]T_0,T_1[$ of $\L u=0$ it holds
 $$ u(x,y,t) \leq \left( \tfrac{t-T_0}{t_0-T_0} \right)^\beta M^{1+
\frac{\Phi(\o)}{h}+\frac{4(t_0-t)}{\theta^2}} 
u(x_0,y_0,t_0),
 $$
where
\begin{equation}
 \Phi(\o) = \int_t^{t_0} \o^2 (s) \, d s.
\end{equation}
\end{proposition}

\medskip

\noindent{\sc Proof.} If $\omega\in L^1([t,t_0]) \setminus L^2([t,t_0])$, then our claim reads as $u(x,y,t)\leq
+\infty$, that is clearly true. We now assume $\omega\in L^2([t,t_0])$. The proof of the proposition is
based on the construction of a Harnack chain, by applying several times Corollary \ref{c1a}. We then first
fix $\theta\in]0,1[$ as in Corollary \ref{c1a}, and we also fix the constant $h=4 \log^2(3/2)$.

{\bf Step 1.} We fix three restrictive assumptions:
\begin{itemize}
\item it holds $t_0-T_0\leq \frac14$ ;
\item the path $\gamma$ is defined on the time interval $[0,t_0-t]$ with  $t_0-t\leq \theta^2(t_0-T_0)$;
\item the function $\Phi(\omega)$ satisfies $\Phi(\omega)\leq h$.
\end{itemize}

We first claim that, under such hypotheses, it holds
\begin{equation}\label{e3a}
 \g(t+s) \in \P_{r}(x_0,y_0,t_0) \quad \text{for every} \
 s \in [0,t_0-t],
\end{equation}
with $r:=\sqrt{t_0-T_0}\leq\frac12$. Indeed, H\"older inequality implies
\begin{equation*}
    \left| \int_t^{t+s} \o(\t) d \t \right| \le \sqrt{s} \left(\int_t^{t+s}
    \o^2(\t) d \t \right)^\frac12 \le \sqrt{h} \sqrt{s} \leq \log (1 +
    \sqrt{s}),
\end{equation*}
for every $s \in [0,t_0-t]\subset[0,\frac14]$. The last inequality follows from concavity of $\log(1+a)$, that implies
$\log(1+a)\geq 2\log(3/2)a$ for $a\in [0,1/2]$ and from the definition of $h$.
We then find
\begin{equation*}
    \left| e^{\int_t^{t+s} \o(\t) d \t} -1 \right| \le  e^{\left|\int_t^{t+s}
    \o(\t) d \t\right|} -1 \le \sqrt{s}
\end{equation*}
for every $s \in [0,t_0-t]$. Thus, by integrating the system \eqref{paramet-gamma}, we obtain
\begin{equation*}
    |x(s) - x_0| \le \sqrt{s} x_0, \quad \text{and} \quad |y(s) - y_0 - s x_0| \le \tfrac{2}{3}s^\frac32 x_0<s^\frac32
x_0
\end{equation*}
for every $s \in [0,t_0-t]$, and \eqref{e3a} is proven. Since $H_r(x_0,y_0,t_0)\subset \R ^+\times \R \times ]T_0,T_1[$
for the definition of $r$, then Corollary \ref{c1a} can be applied, and it holds
$u(x,y,t)\leq M u(x_0,y_0,t_0)$ with $M$ given in Proposition \ref{p1a}.

{\bf Step 2.} We now remove the three hypotheses of Step 1 and prove the main statement. Consider any control
$\omega\in L^2([t,t_0])$ and the corresponding curve $\gamma(.)$. Define the sequence of times
$t<t_k<t_{k-1}<\ldots<t_2<t_1<t_0$ recursively starting from $t_0$ as follows
\begin{equation} \label{eq-recursive}
t_{j+1}=\max\left\{t, \ t_j-\theta^2/4, t_j-\theta^2(t_j-T_0), \   \inf \left\{ s
\mbox{~~s.t.~~}\int_s^{t_j}|\omega(\tau)|^2\,d\tau\leq h\right\}\right\}.
\end{equation}
It is easy to prove that such sequence terminates in a finite number of steps, when the lower boundary $t$ is reached. For simplicity of notation, we denote
$t_{k+1}=t$.

We now define $r_j=\sqrt{t_j-t_{j+1}}/\theta$ , then we note that $r_j \le 1/2$  and
$$H_{r_{j}}(x(t_0-t_{j}),y(t_0-t_{j}),t_{j}) \subset \R ^+\times \R  \times[T_0,T_1],$$
by \eqref{eq-recursive}.
 Moreover, we clearly have $t_j-t_{j+1}\le \theta^2 r_j^2$. By applying Step 1 on the $k+1$ intervals  $[t_{j+1},t_j]$,
it holds
$$
  u(x,y,t)\leq M^{1+k}u(x_0,y_0,t_0).
$$
We point out that the points $(x(t_j),y(t_j), t_j), \ j=1,\ldots k+1$, selected on the path $\gamma(.)$, form a Harnack
chain. Since \eqref{eq-recursive} implies
$$
  k \leq \frac{\int_t^{t_0} |\o(\t)|^2 d \t}{h}+ 4\frac{t_0-t}{\theta^2}+\frac{1}{|\log(1-\theta^2)|}\log\big(
\tfrac{t-T_0}{t_0-T_0}\big),
$$
this concludes the proof of Proposition \ref{p2a}, by setting $\beta := \frac{\log(M)}{|\log(1-\theta^2)|}$.  \hfill
$\square$

\medskip

\begin{remark} \label{rem-L-L0}
Even if $\L$ does not write in the form \eqref{e-Hormanderopeerator}, the lower bound in Proposition \ref{p2a} basically
depends on $\gamma$, that in turns depends on the vector fields $X$ and $Y$ that define $\LY$. This feature depends on
the fact that $\gamma$ is contained in the set $\P_r(z_0)$, where the Harnack inequality holds for both operators $\LY$
and $\L$.
\end{remark}
\medskip

\section{Elements of Stochastic theory}
\label{sec2}
\setcounter{section}{3} \setcounter{equation}{0} \setcounter{theorem}{0}

This section contains some known results about the theory of diffusion processes we need in this work.
We refer to the monograph of Nualart \cite{Nualart}, and Bally \cite{Bally} for an exhaustive presentation of the topic.

Throughout this section, we denote by $C_{l,b}^{\infty}(\R^N,\R)$ the space of smooth functions with
bounded derivatives of any order. Note that the boundedness of the functions is not required. We denote
by $C_p^\infty(\R^N)$ the set of smooth functions $f:\R^N\rightarrow \R$ such that $f$ and all its partial
derivatives have polynomial growth.

We consider the $N$-dimensional Markovian diffusion process $(X_t)_t$ solution of the SDE:
\begin{equation}\label{diffusion}
 dX^i_t=\sum_{j=1}^d\sigma^i_j(X_t)dW^j_t+ F^i(X_t)dt, \qquad i= 1, \dots, N, \quad t \geq 0
\end{equation}
where $W_t=(W_t^1,\ldots, W_t^d)$ is a $d$-dimensional Brownian motion,  $(X_t)_{t\geq 0}$ is a stochastic process on a
probability space $(\Omega,\F,\mathbb{P})$ endowed with the filtration $(\F_t)_{t\geq 0}$ generated by $ (W_t)_{t\geq
0}$ and belonging to the space $L^2([0,\infty) \times \Omega;\mathcal{B}_+\times \F; \lambda \times \mathbb{P})$, where
$\lambda$ stands for the Lebesgue measure in $\R^N$ and $\mathcal{B}_+$ is the Borel $\sigma$-algebra. We assume that
$$F^i, \ \sigma_j^i \ \in C_{l,b}^{\infty}(\R^N,\R)\quad i=1,...,N; \quad j=1,...,d.$$
We denote by $X^x_t$ the solution of the SDE (\ref{diffusion}) with initial condition $X_0^x=x \in \mathbb{R}^N$.

\smallskip

By using the Feynman-Kac representation formula (see for instance Pascucci \cite[chap.9]{Pascucci-book}), one can
state that  the transition density (whenever it exists) $p(x_0,t_0,x,t)$ of the $N$-dimensional process
\eqref{diffusion} satisfies the Fokker-Planck equation:
\begin{equation}\label{FP}
  \sum_{i,j=1}^da_{ij}(x)\partial_{x_ix_j}p(x_0,t_0,x,t)+
\sum_{i=1}^nF_i(x)\partial_{x_i}p(x_0,t_0,x,t)+\partial_tp(x_0,t_0,x,t)=0
\end{equation}
where
$$ a_{ij}(x)=\frac{1}{2}\sum_{k=1}^d \sigma_{ik}(x)\sigma_{jk}(x).$$
Specifically, the function
\begin{equation}\label{link}
  u(x,t)=\mathbb{E}[\phi(X_T)|X_t=x]=\int_{\R^N}\phi(\xi) p(\xi,T;x,t)d\xi
\end{equation}
is a solution of the Cauchy problem for \eqref{FP} with prescribed bounded continuous final condition $\varphi$.
Moreover, $p$ satisfies the identity
\begin{equation}\label{reproduction_p}
  p(x_0,t_0;x,t)=\int_{\R^N}p(x_0,t_0;\xi,\tau)p(\xi,\tau;x,t)d \xi, \qquad t<\tau<t_0.
\end{equation}

In the sequel of this section we recall the results of the Stochastic Theory which guarantee the existence of the
transition density $p(x_0,t_0,x,t)$.

\subsection{Elements of Malliavin Calculus}
We consider the space of functions $\mathcal{H}=L^2([0,T], \R^d)$. For each $h(t)=(h^1(t),...,h^d(t)), \in \mathcal{H}$
we introduce the
Gaussian random variable:
$$W(h)= \sum_{j=1}^d \int_0^T h^j(t)dW_t^j. $$
We denote by $\mathcal{S}$ the class of $n$-dimensional simple functions of Brownian motion of the form:
$$ F=f(W(h_1),...,W(h_n)), \quad f \in C_p^\infty(\R^n, \R), \quad h_1,...,h_n \in \mathcal{H}.$$
For every $F \in \mathcal{S}$ we define the \emph{Malliavin derivative} $(D_tF)_{t \in [0,T]}$ of $F$ as the
$\R^d$-dimensional (non adapted) process:
$$D_tF=\sum_{i=1}^n\frac{\partial f}{\partial x_i}(W(h_1),...,W(h_n)))h_i(t).$$
Each $h_i(t)=(h_i^1(t),...,h_i^d(t))$ has $d$ components and we write $D_t^jF$ for the $j^{th}$ component of $D_tF$,
$j=1,...,d$.
We introduce the Sobolev norm:
\begin{equation}\label{normsob}
\|F\|_{1,p}=\Big[\mathbb{E}(|F|^p)+\mathbb{E}\big(\big|DF\big|^p\big)\Big]^{1/p}
\qquad \mbox{where}\qquad
\big|DF\big|=\left(\int_0^T|D_tF|^2dt\right)^{1/2}.
\end{equation}
It is possible to show that
the operator $D: \mathcal{S} \rightarrow L^p(\Omega,L^2[0,T])$ is closable with respect to the norm $\| \cdot\|_{1,p}$.
We denote by $\mathbb{D}^{1,p}=Dom(D)$ its domain, which is the completion of $\mathcal{S}$ with respect to the norm
$\| \cdot\|_{1,p}$.

\smallskip

 Let $\alpha=(j_1,...,j_k)$ be a multi-index of length $k$,  we define  the $k^{th}$-order derivative as the random
vector on $[0,T]^k \times \Omega$ with coordinates:
$$D^\alpha_{t_1,...,t_k}F=D_{t_k}^{j_k} \cdots D_{t_1}^{j_1}F.$$
We introduce the Sobolev norm:
\begin{equation}\label{normsob1}
\|F\|_{k,p}=\Big[\mathbb{E}(|F|^p)+\sum_{j=1}^k\mathbb{E}(|D^{(j)}F|^p)\Big]^{1/p}
\end{equation}
where
$$ |D^{(j)}F|=\sum_{|\alpha|=j}\left(\int_{[0,T]^k}|D_{t_1,...,t_k}^{\alpha}F|^2dt_1\ldots dt_k\right)^{1/2}$$
We denote by $\mathbb{D}^{k,p}$ the completion of $\mathcal{S}$ with respect to the norm $\|\cdot\|_{k,p}$ and finally
we denote by
$$\mathbb{D}^\infty= \bigcap_{k,p \geq 1}\mathbb{D}^{k,p}.$$
We introduce now  the Malliavin covariance matrix of the random vector $F=(F^1,...,F^N)$ derivable in Malliavin sense.
\begin{definition}
Let $ F=(F^1, \ldots , F^N)$ be a random vector which is derivable in Malliavin sense.\\ We define the \emph{Malliavin
Covariance Matrix} of the random variable $F$ as follows:
\begin{equation}\label{MallMatrix}
  \gamma_{F}^{ij}=\langle DF^i, DF^j\rangle=\sum_{k=1}^d \int_0^TD_s^k F^i \times D_s^kF^j ds \quad i,j=1,...,N
\end{equation}
We say that $F$ is \emph{non-degenerate} if its Malliavin covariance matrix satisfies
\begin{equation}\label{non-deg}
\mathbb{E}(|\text{\rm{det}} \gamma_{F}|^{-p})<\infty, \quad \forall p \in \mathbb{N},\ \ \forall t>0
\end{equation}
\end{definition}
The \emph{non-degeneracy} \eqref{non-deg} condition is necessary to ensure that the law of the random vector $F$ exists
and is
absolutely continuous with respect to the Lebesgue measure. We refer to \cite{Nualart}, Chapter 2, for the following
proposition
\begin{proposition}[Hirch-Bouleau]\label{Hirch-Boleau}
Let $t \in [0, +\infty)$ be fixed and let $X_t=(X_t^1,...,X_t^n)$ a random variable satisfying (\ref{diffusion}).
If each $X_t^i \in \mathbb{D}_{\loc}^{1,p}$ with $p>1$ and if $\g_{X_t}$ satisfies the \emph{non degeneracy condition}
(\ref{non-deg}) almost surely, then the law of $X_t$ is absolutely continuous with respect to the Lebesgue measure on
$\R^N$, that is
$$P_{X_t}(dx)=p_{X_t}(x)dx.$$
\end{proposition}
\subsection{Malliavin Theorem and H\"{o}rmander condition}
In this section we recall the Malliavin Theorem for a diffusion process \eqref{diffusion}.

\begin{theorem}[Malliavin \cite{Malliavin}]\label{Malliavin}
Consider the $n$-dimensional diffusion process (\ref{diffusion}) and suppose that $F^i, \ \sigma_j^i \  \in
C_{l,b}^{\infty}$.
\begin{description}
\item[i)] Then for every $t>0$, $X_t$ belongs to $\mathbb{D}^\infty$ and
\begin{equation}\label{M1}
 \| X_t^x\|_{k,p}\leq c_{k,p}(t)(1+|x|)^{\beta_{k,p}}
\end{equation}
where $\beta_{k,p} \in \mathbb{N}$ and $c_{k,p}(t)$ is a constant which depends on $k,p, t$ and on the bounds of the
derivatives of $b, \sigma$ up to order $k$.
\item[ii)] Suppose that H\"{o}rmander condition \eqref{e-Hormandercondition} holds true.  Then there exist a function
$C_{k,p}(t)$ and some constants $n_k, m_k \in \mathbb{N}$ such that the non-degeneracy condition \eqref{non-deg} is
satisfied. Moreover
\begin{equation}\label{M2}
  \|
  \, (\gamma_{X_t^x})^{-1}\|_p \leq \frac{C_{k,p}(t)(1+|x|)^{m_k}}{t^{n_k/2}}.
\end{equation}
The function $t\rightarrow C_{k,p}(t)$ is increasing. In particular, the right hand side in \eqref{M2} blows up as
$t^{-n_k/2}$ as $t
\rightarrow 0$.
\item[iii)] Suppose that the H\"{o}rmander condition \eqref{e-Hormandercondition} holds true and $F^i, \ \sigma_j^i \
\in C_{l,b}^{\infty}$.
Then for every $t>0$ the law of $X_t^x$ is absolutely continuous with respect to the Lebesgue measure and the
transition density $ y \mapsto p(y,t;x,t_0)$ is a $C^\infty$ function. Moreover, if $b, \sigma$ are bounded, one has
\begin{gather}\label{M3}
  p(y,t;x,t_0) \leq \frac{C_0(1+|x|)^{m_0}}{t^{n_0/2}}\exp \left(- \frac{D_0(t)|y-x|^2}{t}\right)\\
  |D_y^\alpha p(y,t;x,t_0)| \leq \frac{C_\alpha(1+|x|)^{m_\alpha}}{t^{n_\alpha/2}}\exp
\left(-\frac{D_\alpha(t)|y-x|^2}{t}\right)
\end{gather}
where all above constants depend on the step for which H\"{o}rmander condition holds true and the functions
$C_0,D_0,C_\alpha,D_\alpha$ are increasing functions of $t$.
\end{description}
\end{theorem}

We now consider the operator $\L$ in \eqref{e-main}, assuming that the coefficients $a,b$ only depend on $x,y$ and are
bounded $C^{\infty}(\R^2)$ functions. We denote by
\begin{equation}\label{Lnondiv}
  \mathcal{L}=a(x,y)x^2\p_{xx}+\big(a_x(x,y)x+a(x,y)+b(x,y)\big)x\p_x+x\p_y.
\end{equation}
and from \eqref{FP} we have that $\mathcal{L}+\p_t$ is the infinitesimal generator of the process
\begin{equation}\label{process_main}
  \left\{
    \begin{array}{ll}
      dX_t=\mu(X_t,A_t)X_tdt+ \sigma(X_t,A_t)X_tdW_t \\
      dY_t=X_tdt.
    \end{array}
  \right.
\end{equation}
with
$$a(x,y)=\frac{\sigma^2(x,y)}{2}, \qquad b(x,y)+\frac{\sigma^2(x,y)}{2}+\sigma(x,y)\sigma_x(x,y)x=\mu(x,y).$$
It is simple to show that the process $(X_t,Y_t)_{t \ge 0}$ belongs to the space $C_{l,b}^\infty(\R^+\times\R)$,
provided that $\p_x(xa(x,y))$ is bounded. Moreover, the operator \eqref{Lnondiv} satisfies the H\"{o}rmander Condition,
then the
density $p$ of the process  $(X_t,Y_t)_{t \ge 0}$ exists in view of \textbf{i)} and \textbf{ii)} of Theorem
\ref{Malliavin} and Proposition \ref{Hirch-Boleau}. Point \textbf{iii)} of Theorem \ref{Malliavin} yields the
smoothness of $p$.

\medskip

The following proposition  summarizes the results about the fundamental solution of $\L$ we have obtained in this
Section.

\begin{proposition} \label{prop-ex-fsmc}
Let $a = a(x,y), b = b(x,y) \in C^\infty(\R^+ \times \R)$, with $a, b$ and $\partial_x(x a(x,y,t))$ bounded. Suppose
that $\inf a > 0$. Then, there exists a smooth fundamental solution of $\L$. Moreover, for every $(x,y,t),
(x_0,y_0,t_0), (\xi,\eta,\tau)$ belonging to $\R^+\times\R^2$ with $t>\tau>t_0$, it holds the following properties
\begin{enumerate}
  \item {Support of} $\G$:
  \begin{equation}\label{support}
  \Gamma(x,y,t;\xi,\eta,\tau)=0 \qquad \mbox{whenever} \qquad t \le \tau \ \ \mbox{or} \ \ y \ge \eta;
\end{equation}
  \item {Reproduction property}:
  \begin{equation}\label{reproduction}
  \G(x,y,t;x_0,y_0,t_0)= \int_{\R^+ \times \R}\G(x,y,t;\xi,\eta,\tau)\G(\xi,\eta,\tau;x_0,y_0,t_0)d\xi d \eta,
\end{equation}
  \item  Integrals of $\Gamma$: \begin{equation}\label{integrals}
   \int_{\R^+ \times \R}\G(x,y,t;\xi,\eta,\tau)d\xi d \eta=1, \qquad \int_{\R^+ \times \R}\G(x,y,t;\xi,\eta,\tau)dx dy= \bar{C};
\end{equation}
where $\bar{C}$ is a positive constant depending on $t-\tau$ and $\bar{C} \rightarrow 1$ when $t\rightarrow \tau$.
\end{enumerate}
\end{proposition}

\medskip

\noindent {\sc Proof.}
 Malliavin Calculus provides us with the existence of a smooth probability density $p(x_0,y_0,t_0;x,y,t)$ for the
process \eqref{process_main}. By setting
\begin{equation}\label{link_Fund_prob}
\Gamma(x,y,t;\xi,\eta,\tau)=p(\xi,\eta,T-\tau,x,y,T-t).
\end{equation}
it is easy to check that \eqref{link_Fund_prob} defines a smooth Fundamental solution for $\L$ in the sense of the
Definition \ref{def_fundsol}.
The relation \eqref{support} simply follows from \eqref{process_main}, as the process $(X_t)_{t \ge 0}$ is
positive. The reproduction property
\eqref{reproduction} follows from \eqref{reproduction_p} and \eqref{process_main}. Moreover, the first property \eqref{integrals} follow from \eqref{process_main} and by the fact that $p$ is the transition probability density of a Markovian process.\\
In order to prove the second property of \eqref{integrals}, we consider the adjoint operator $\L^*$ of $\L$ as in \eqref{main_adj}. Let rewrite $\L^*$ in the following form
\begin{align*}\label{divergence_form1}
\L^*=a(x,y)x^2\p_{xx}+\big(a_x(x,y)x+&a(x,y)-b(x,y)\big)x\p_x-x\p_y+\\
&-(b(x,y)-a(x,y)-xa_x(x,y)-b_x(x,y))+\p_t.
\end{align*}
whose fundamental solution is $\G^* (\xi,\eta,\tau;x,y,t)$ with $t> \tau$.\\
Denote by $\tilde{\Gamma}^*(\xi,\eta, \tau;x,y,t)$ the fundamental solution of the operator
\begin{equation}\label{}
\tilde{\L}^*=a(x,y)x^2\p_{xx}+\big(a_x(x,y)x+a(x,y)-b(x,y)\big)x\p_x-x\p_y+\p_t.
\end{equation}
Note that $\tilde{\Gamma}^*(\xi,\eta, \tau;x,y,t)$ agrees with the probability density of the Stochastic Process:
\begin{equation*}\label{process_main_adj}
  \left\{
    \begin{array}{ll}
      dX_t=\mu^*(X_t,A_t)X_tdt+ \sigma(X_t,A_t)X_tdW_t, & X_\tau=\xi \\
      dY_t=-X_tdt, & Y_\tau=\eta.
    \end{array}
  \right.
\end{equation*}
with
$$a(x,y)=\frac{\sigma^2(x,y)}{2}, \qquad -b(x,y)+\frac{\sigma^2(x,y)}{2}+\sigma(x,y)\sigma_x(x,y)x=\mu^*(x,y).$$
Therefore, it follows that:
\begin{equation}\label{integral2}
  \int_{\R^+ \times \R}\tilde{\Gamma}^*(\xi,\eta,\tau;x,y,t)dx dy =1
\end{equation}
Let $C_2$ denote the maximum of the function $|b(x,y)-a(x,y)-xa_x(x,y)-b_x(x,y)|$ and  consider the following operators
\begin{gather}
\L^*_1=a(x,y)x^2\p_{xx}+\big(a_x(x,y)x+a(x,y)-b(x,y)\big)x\p_x-x\p_y+C_2+\p_t\\
\L^*_2=a(x,y)x^2\p_{xx}+\big(a_x(x,y)x+a(x,y)-b(x,y)\big)x\p_x-x\p_y-C_2+\p_t.
\end{gather}
Observe that the functions ${\G}^*_1(\xi,\eta, \tau;x,y,t)=e^{-C_2(t-\tau)}\tilde{\G}^*(\xi,\eta, \tau;x,y,t)$ and ${\G}^*_2(\xi,\eta, \tau;x,y,t)=e^{C_2(t-\tau)}\tilde{\G}^*(\xi,\eta, \tau;x,y,t)$ are the fundamental solution of $\L^*_1 v=0$ and $\L^*_2 v=0$, respectively. Moreover, for every non negative function $g(x,y)$  continuous and bounded on $\R^+\times \R$, if we consider
\begin{align}
u_1(\xi,\eta,\tau)=\int_{\R\times \R^2} \G^*_1 (\xi,\eta,\tau;x,y,t)g(x,y) d x d y, \nonumber \\
u(\xi,\eta,\tau)=\int_{\R\times \R^2} \G^*(\xi,\eta,\tau;x,y,t)g(x,y) d x d y, \nonumber \\
u_2(\xi,\eta,\tau)=\int_{\R\times \R^2} \G^*_2(\xi,\eta,\tau;x,y,t)g(x,y) d x d y, \nonumber
\end{align}
it holds that
\begin{equation}\label{soprasol_sottsol}
\L^* u_1 (\xi,\eta,\tau) \leq 0, \quad \L^* u (\xi,\eta,\tau) = 0, \quad \L^* u_2 (\xi,\eta,\tau) \geq 0,
\end{equation}
 In view of \eqref{soprasol_sottsol} and by using the comparison principle, we obtain:
\begin{equation}\label{soprasol_sottsol2}
 \G^*_1 (\xi,\eta,\tau;x,y,t) \leq \G^* (\xi,\eta,\tau;x,y,t)\leq \G^*_2 (\xi,\eta,\tau;x,y,t).
\end{equation}
The assertion simply follows from the fact that, in view of \eqref{integral2}, we have
\begin{eqnarray}\label{integral3}
  \int_{\R^+ \times \R}{\Gamma}^*_1(\xi,\eta,\tau;x,y,t)dx d y=e^{-C_2(t-\tau)}, \quad \nonumber
  \int_{\R^+ \times \R}{\Gamma}^*_2(\xi,\eta,\tau;x,y,t)dx d y=e^{C_2(t-\tau)}.
\end{eqnarray}
and $\G^* (\xi,\eta,\tau;x,y,t)=\G (x,y,t; \xi,\eta,\tau)$.
\hfill $\square$.

\section{The Optimal Control Problem and The Lower Bound}\label{sec4}
\setcounter{section}{4} \setcounter{equation}{0} \setcounter{theorem}{0}

In this section we formulate the control problem suitable to find the optimal lower bound for the positive solutions of
$\L u = 0$. With this aim, we recall that $\LY$ can be written in the form \eqref{e-Hormanderopeerators}, as $\LY = X^2
+Y$, where $X$ and $Y$ are defined in \eqref{e-vectorfields}.

\subsection{The Pontryagin Maximum Principle}
In this section, we recall the Pontryagin Maximum Principle \cite{PMP}. We will then apply it to the optimal control
problem \eqref{contr-ott-0} in Section \ref{s42}, and it will give us optimal lower bounds for the positive solutions
of
$\L u = 0$. We use here the notations in the general setting suitable for the study of operators $\widetilde \L$
defined
in \eqref{e-Hormanderopeerators} that include, as a particular case, the one studied in this work.

In this section the time variable $t$ in $(x,t)\in\R^N\times [0,+\infty)$ is dropped. Let then $\Omega \subset \R^N$ be
an open set, $F_0, F_1, \dots, F_m: \Omega \to \R^N$ be smooth vector fields, and the final time $T$ be fixed. We
consider the following {\it optimal control problem}:
\begin{eqnarray}
\dot q&=&F_0(q)+\sum_{i=1}^m \omega_i F_i(q)\,,~~\o_i\in\R\,,
~~\int_0^T
\sum_{i=1}^m \omega_i^2(t)~dt\to\min,~~q(0)=q_0,~~q(T)=q_1.
\label{sopra}
\end{eqnarray}

For such optimal control problem, the Pontryagin Maximum Principle provides a first-order condition for the minimizing
controls $\omega(.)$ and the corresponding trajectories $q(.)$. We now recall its statement in the particular case in
which variables and controls belong to the Euclidean spaces $\R^n,\R^m$, respectively. For a more general statement on
manifolds, see e.g. \cite{agra-book}.

\begin{theorem}[PMP for the problem \eqref{sopra}] \label{t-PMP}
\it Consider the minimization problem {\rm \eqref{sopra}}, in the class of Lipschitz continuous curves,
where  $F_i$, $i=0,\ldots, m$ are smooth vector fields on $\R^N$ and the final time $T$ is fixed.
Consider  the map $H:\R^N\times\R^N \times\R\times \R^m\to \R$ defined  by
\begin{eqnarray}
H(q,\lambda,p_0,\omega)&:=&\left\langle\lambda, F_0+\sum_{i=1}^m \omega_i F_i(q)\right\rangle+p_0 \sum_{i=1}^m
\omega_i^2.
\end{eqnarray}
If the curve $q(.):[0,T]\to \R^N$ corresponding to the control $\omega(.):[0,T]\to\R^m$
is optimal, then there exist a Lipschitz continuous {covector}
$\lambda(.):s\in[0,T]\mapsto \lambda(s)\in \R^N$ and a constant $p_0\leq 0$ such that:
\begin{itemize}
\item the pair $(\lambda(s),p_0)$ is never vanishing;
\item the optimal control $\omega(s)$ satisfies
$$H(q(s),\lambda(s),p_0,\omega(s))=\max_{\nu\in\R^m}H(q(s),\lambda(s),p_0,\nu);$$
\item  for a.e. $s\in [0,T]$ it holds
\begin{equation}
\begin{cases}
\dot q(s)={\frac{\partial H
}{\partial \lambda}(q(s),\lambda(s),p_0,\omega(s))},&\\
\dot \lambda(s)=-{\frac{\partial H
}{\partial q}(q(s),\lambda(s),p_0,\omega(s))}.
\end{cases}
\label{ePMP}
\end{equation}
\end{itemize}
The Hamiltonian $H^*(q,\lambda,p_0):=\max_{\nu\in\R^m} H(q,\lambda,p_0,\nu)$ is called the {\bf maximized Hamiltonian}.

Solutions to the system \eqref{ePMP} are called {\bf extremals}. When $p_0=0$, they are called {\bf abnormal extremals},
while when $p_0<0$ they are called {\bf normal extremals}.
\end{theorem}

\begin{remark} The original statement \cite{PMP} of the Pontryagin Maximum Principle provides optimal controls in the
space $L^\infty([0,T],\R^m)$. Instead, we are interested in optimal controls in the larger space $L^1([0,T],\R^m)$. For
this reason, we aim to apply a generalized version of the Pontryagin Maximum  Principle, such as the one stated in
\cite[Chap. 6]{vinter}. For our optimal control problem, such generalized version has a statement completely equivalent
to Theorem \ref{t-PMP}.
\end{remark}

\subsection{Application of the Pontryagin Maximum Principle to the problem \eqref{contr-ott-0}}\label{s42}
In this section we apply the Pontryagin Maximum Principle to our problem \eqref{contr-ott-0}. Note that the
terminal point of the $\L$-admissible path considered in \eqref{contr-ott-0} is $(1,0,0)$, we give here the formulation
for any end-point $(x_0,y_0, t_0) \in \R^+ \times \R^2$. In accordance with the notation used for the fundamental
solution of $\L$, we denote the starting point of the path by $(x_1,y_1,t_1) \in \R^+ \times \R^2$, with $t_1 > t_0$.
\begin{eqnarray} \label{prob42a}
&  \left\{
    \begin{array}{ll}
      \dot{x}(s)=\o(s)x(s) & \\
      \dot{y}(s)=x(s) &\hbox{$ 0 \le s \le T$,}\\
      \dot t(s)=-1,&
    \end{array}
  \right.
  \\
&{(x,y,t)(0)=(x_1,y_1,t_1),\qquad (x,y,t)(T)=(x_0,y_0,t_0).} \nonumber
\end{eqnarray}
We first observe that such optimal control problem is invariant on the Lie group $\R^+ \times \R^2$ endowed with the
operation \eqref{sx}. We recall that optimal control problems on Lie group with invariant vector fields satisfy useful
invariance properties, that permit to have simpler solutions of the Pontryagin Maximum Principle, eventually leading to
complete synthesis for specific problems, see e.g. \cite{boscain2008invariant}.
In our specific problem, it is sufficient to observe the following invariance property for the solution of
\eqref{prob42a}. Consider a control $\omega( \, \cdot \, )$ steering $(x_1,y_1,t_1)$ to $(x_0,y_0,t_0)$ with the
trajectory $(x(s),y(s),t(s))$. Then the same control $\omega(.)$ steers $(x_0,y_0,t_0)^{-1}\circ (x_1,y_1,t_1)$ to
$(1,0,0)$. This can be proved by observing that the trajectory $(x_0,y_0,t_0)^{-1} \circ (x(s),y(s),t(s))$ is a
solution of \eqref{prob42a} with the same control $\omega(.)$. Since the cost depends on the control only, then the two
trajectories have the same cost, hence
\begin{equation}\label{transl-phi}
{\Psi(x_1,y_1,t_1;x_0,y_0,t_0)=\Psi((x_0,y_0,t_0)^{-1}\circ(x_1,y_1,t_1);1,0,0)}
\end{equation}
As a consequence, we will now fix the final condition $(x_0,y_0,t_0)=(1,0,0)$ in the optimal control problem
\eqref{prob42a}, then using the invariance property to solve it with a general initial condition.

The constraint $\dot t=-1$ implies that $\L$-admissible paths satisfy $t(s)=t_1-s$, hence $T= t_1 - t_0$. Then, in
the sequel we drop the time variable, we set $T := t_1 - t_0$, and we denote
\begin{equation} \label{cost-phi}
  \Psi(x_1, y_1, t_1; x_0, y_0, t_0) = \inf_{\omega \in L^1([0,t_1 - t_0])} \int_0^{t_1 - t_0}\omega^2(\tau) d\tau,
\end{equation}
where $\omega \in L^1([0, t_1 - t_0])$ is such that \eqref{prob42a} holds true.

For the above reasons, the optimal control problem \eqref{prob42a}, \eqref{cost-phi} now reads as follows:
\begin{eqnarray}\label{prob-ott}
\Psi(x_1, y_1, t_1; 1,0,0) = \min_{\omega \in L^1([0,t_1])} \int_0^{t_1}\omega^2(\tau) d\tau
\quad \mbox{subject to constraint}
\\ \nonumber
\\ \nonumber \quad \left\{
             \begin{array}{ll}
               \dot x(s)=\omega(s) x(s), & \hbox{$x(0)= x_1, \quad x(t_1)=1$,} \\
               \dot y(s)=x(s), & \hbox{$y(0)= y_1, \quad y(t_1)=0$.}
             \end{array}
           \right.
\end{eqnarray}
To simplify the notation, in the sequel we agree to set ${\Psi(x_1, y_1, t_1) := \Psi(x_1, y_1, t_1; 1,0,0)}$.

We now solve such problem. As a by-product, we show that we can always steer $(x_1,y_1)$ to $(x_0,y_0)$ in time $T$,
when $y_1 < y_0$. This implies that there exists a control $\omega$ steering $(x_1,y_1,t_1)$ to $\gamma(t_1-t_0) =
(x_0,y_0,t_0)$, as we stated in the proof of Proposition \ref{supp-density}.

We now apply the Pontryagin Maximum Principle to problem \eqref{prob-ott}. The Hamiltonian of the problem
\eqref{prob-ott} is
\begin{equation}\label{ham}
  H(x,y,\lambda_1,\lambda_2,p_0,\omega) = \lambda_1 x \omega + \lambda_2 x + p_0 \omega^2,
\end{equation}
where $(\lambda_1,\lambda_2)$ are the coordinates of the covector $\lambda$.

We first remark that Problem \eqref{prob-ott} admits no abnormal extremals. Indeed, assume by contradiction $p_0=0$ in
\eqref{ham}. Then
$$
  H(x,y, \lambda_1,\lambda_2,p_0,\omega) = \lambda_1 x \omega + \lambda_2 x
$$
Recall that $x > 0$. Hence, the maximization of the Hamiltonian is equivalent to
$$
  \frac{\partial H}{\partial\omega} (x,y,\lambda_1,\lambda_2,p_0,\omega) = 0
  \quad \Rightarrow  \quad \lambda_1(s)=0, \quad \forall s \in [0,t_1].
$$
Moreover, using the fact that $\lambda_1(s)=0$ for all $s\in[0,t_1]$, it holds
$$
  \dot{\lambda}_1 (s) = -\frac{\partial H}{\partial x} (x,y,\lambda_1,\lambda_2,p_0,\omega)
=-\lambda_1(s)\omega(s)-\lambda_2(s)=0,
$$
hence $\lambda_2(s)=0$ for every $s \in [0,t_1]$.
We conclude that
$$
  (\lambda_1(s), \lambda_2(s),p_0)=(0,0,0) \quad \text{for every} \quad s \in [0,t_1].
$$
This is in contradiction with the fact that $ (\lambda_1(s), \lambda_2(s),p_0)$ is always non-vanishing.

Since no abnormal extremals occur, we choose $p_0=-\frac12$. We then compute the optimal control as the unique
minimizer of $H\left(x,y,\lambda_1,\lambda_2,-\frac12,\omega\right)$, that is
\begin{equation}\label{probl-omega}
\omega(s)=\lambda_1(s) x(s),
\end{equation}
and the maximized Hamiltonian is
\begin{equation}\label{prob-H}
H^*(x,y, \lambda_1,\lambda_2,p_0)=\frac{1}{2} \lambda_1^2x^2 +\lambda_2x.
\end{equation}
The corresponding Hamiltonian system reads as
\begin{equation}\label{probl-2}
\left\{
  \begin{array}{ll}
   \dot{x}(s)=\lambda_1(s) x^2(s) \\
     \dot{y}(s)=x(s)\\
    \dot{\lambda}_1(s)=-\lambda_1^2(s)x(s)-\lambda_2(s) \\
    \dot{\lambda}_2(s)=0
  \end{array}
\right.
\end{equation}
In the sequel, we choose the parameters
$$
 { k := \lambda_1(t_1)\qquad\mbox{~~ and ~~}\qquad c := \lambda_2(t_1)}
$$
as the {final} condition for each extremal, that is uniquely determined by being the solution of \eqref{probl-2} with
final condition $(x,y,\lambda_1,\lambda_2)(t_1)=(1,0,k,c)$. Note that, by the last equation in \eqref{probl-2}, we have
$\lambda_2(s) = c$ for every $s \in [0, t_1]$. Furthermore, the value of the Hamiltonian is a constant of motion, fixed
by the final data. From now on, we then fix
\begin{equation}\label{integrprimi}
      E:={\lambda_1^2(s)x^2(s)}+ 2 \lambda_2(s)x(s)=k^2+2c.
\end{equation}

Moreover, by recalling the explicit expression for the optimal control \eqref{probl-omega} and $\dot y=x$, we have the
following expression of the cost for extremals:
\begin{equation}
  C(\omega(\, \cdot \,))=\int_0^{t_1} \omega(s)^2\,ds=\int_0^{t_1} \lambda_1^2(s)x^2(s)\,ds=
  \int_0^{t_1} (E-2c \dot y(s))\,ds = E{t_1} + 2c y_1.
\end{equation}

\medskip

We now describe the explicit solutions to \eqref{probl-2}, as a function of the final value of the Hamiltonian
$E=k^2+2c$. For simplicity, we consider the space variable $(x,y)$ only. We have three cases:
\begin{enumerate}
\item {\bf $E=0$:} it holds $(x(s),y(s)) = \left(\frac{4}{(k(t_1-s)+2)^2},-\frac{2(t_1-s)}{k(t_1-s)+2}\right)$;
\item {\bf $E>0$:} it holds
$$(x(s),y(s)) = \left( \tfrac{E}{\left(\sqrt E \cosh \big(\tfrac{t_1-s}{2}\sqrt E \big)+k
\sinh\big(\tfrac{t_1-s}{2}\sqrt E \big) \right)^2}, \tfrac{-2\sinh\big(\tfrac{t_1-s}{2}\sqrt E \big)}{\sqrt
E\cosh\big(\tfrac{t_1-s}{2} \sqrt E \big)+k\sinh\big(\tfrac{t_1-s}{2} \sqrt E \big)}\right);$$

\item {\bf $E<0$:} it holds
$$\!\!\!(x(s),y(s))  = \! \left(\tfrac{-E}{\left(\sqrt {-E} \cos \big(\tfrac{t_1-s}{2}\sqrt{- E} \big)+
    k \sin\big(\tfrac{t_1-s}{2}\sqrt{- E} \big) \right)^2} , \tfrac{-2\sin\big(\tfrac{t_1-s}{2}\sqrt{- E}
\big)}{\sqrt{- E}\cos \big(\tfrac{t_1-s}{2} \sqrt {-E}\big)+k\sin\big(\tfrac{t_1-s}{2} \sqrt {-E} \big)}\right) \!\!,$$
    where the trajectory is defined on the whole time interval $s\in[0,t_1]$ when $E>-\frac{\pi^2}{T^2}$ only.
\end{enumerate}

The three cases can be unified by using the function $g$ defined in \eqref{func_g} and observing that it always holds
\begin{equation}\label{parab2}
y(s)=-g\left(\frac{E(t_1-s)^2}{4}\right) (t_1-s) \sqrt{x(s)}.
\end{equation}

\medskip

We are now ready to prove the invariance properties of $\Psi$.

\medskip

\noindent{\sc Proof of Proposition \ref{prop-Psi}.}\ The proof of \eqref{tr-inv-psi} is a direct consequence of
\eqref{transl-phi}. In order to prove \eqref{dil-inv-psi} we introduce another symmetry of the problem. Consider an
extremal of \eqref{prob-ott} steering $(x,y)$ to $(1,0)$ in time $t$, with a final covector parametrized by $(k,c)$,
hence with Hamiltonian $E=k^2+2c$ and cost $C = E \, T+2cy_1$. Fix now $r>0$: the extremal ending to $(1,0)$ with final
covector $(r k, r^2 c)$ steers $\left(x,\frac{y}{r}\right)$ to $(1,0)$ in time $\frac{t}{r}$. Moreover, the
Hamiltonian is $r^2 E$ and the cost is $r \, C$. The proof is a direct consequence of the explicit expression of
solutions of \eqref{probl-2}. As a consequence, a trajectory parametrized by $(k,c)$ steering $(x,y)$ to $(1,0)$ in time
$t$ is optimal if and only if the trajectory parametrized by $(r k, r^2 c)$ steering $\left(x,
\frac{y}{r}\right)$ to $(1,0)$ in time $\frac{t}{r}$ is optimal too. Combining this with \eqref{transl-phi}
we get the property
\begin{align*}
\Psi(x_1,y_1,t_1;x_0,y_0,t_0)& = \Psi\big(\tfrac{x_1}{x_0},\tfrac{y_1-y_0}{x_0},t_1-t_0;1,0,0\big) \\ 
  &=\tfrac{1}{r}\Psi\big(\tfrac{x_1}{x_0},\tfrac{y_1-y_0}{r x_0},
\tfrac{t_1-t_0}{r};1,0,0\big)=\tfrac{1}{r}\Psi\big(x_1,\tfrac{y_1}{r},\tfrac{t_1}{r};
  x_0,\tfrac{y_0}{r},\tfrac{t_0}{r}\big) 
\end{align*}
This proves  \eqref{dil-inv-psi}. \hfill $\square$

\medskip

In view of  \eqref{tr-inv-psi} and \eqref{dil-inv-psi}, with no loss of generality, from now on we consider the problem
of steering $(x_1,y_1)$ to $(1,0)$ with fixed final time $t_1=2$. First observe that, since $g$ is a $C^\infty$,
strictly increasing function, from \eqref{parab2} we find the unique value for the prime integral $E$ for which it holds
$(x(0),y(0))=(x_1,y_1)$, that is
\begin{equation}\label{e-E}
E=\frac{4}{t_1^2}g^{-1}\left(-\frac{y_1}{t_1\sqrt{x_1}}\right)=g^{-1}\left(-\frac{y_1}{2\sqrt{x_1}}\right).
\end{equation}

It also clearly gives the basic relation $c=\frac{E-k^2}{2}$, hence $c$ is uniquely determined by $k$. Then, the cost
of
the corresponding extremal is
\begin{equation}\label{cost}
C=2E+y_1(E-k^2)=(2+y_1)E-y_1k^2.
\end{equation}
We now compute the value of $k$ by imposing the initial condition on the second component only, i.e. $y(0)=y_1$. It
holds:
\begin{itemize}
\item for $y_1=-2\sqrt{x_1}$, the unique extremal satisfying $y(0)=y_1$ has final covector $k=-\frac{y_1+2}{y_1}$
and the optimal cost is $C=\frac{(y_1+2)^2}{y_1}$.

\item for $y_1<-2\sqrt{x_1}$, the unique extremal satisfying $y(0)=y_1$ has final covector
$$k=-\sqrt{E}\big(\coth(\sqrt{E})\big)-\tfrac{2}{y_1}=\frac{\sqrt{Ey_1^2+4x_1}-2}{y_1}$$
and the optimal cost is $C=2\frac{Ey_1-2x_1-2+2\sqrt{4x_1+Ey_1^2}}{y_1}.$

\item for $y_1>-2\sqrt{x_1}$, the unique extremal satisfying $y(0)=y_1$ has final covector
$$k=-\sqrt{-E}\big(\cot(\sqrt{-E})\big)-\tfrac{2}{y_1}$$
Since $-\pi^2<E<0$, we find
$$\left\{
    \begin{array}{ll}
      k=\frac{\sqrt{Ey_1^2+4x_1}-2}{y_1}, & \hbox{if $-\pi^2/4\le E<0$;} \\
      k=-\frac{\sqrt{Ey_1^2+4x_1}+2}{y_1}, & \hbox{if $-\pi^2<E<-\pi^2/4$,}
    \end{array}
  \right.
$$
and the expression of the optimal cost is
$$\left\{
    \begin{array}{ll}
     C=2\frac{Ey_1-2x_1-2+2\sqrt{4x_1+Ey_1^2}}{y_1}, & \hbox{if $-\pi^2/4\le E<0$;} \\
     C=2\frac{Ey_1-2x_1-2-2\sqrt{4x_1+Ey_1^2}}{y_1}, & \hbox{if $-\pi^2<E<-\pi^2/4$.}
    \end{array}
  \right.
$$
\end{itemize}
In conclusion,  we have that
 the unique extremal satisfying $y(0)=y_1$ has final covector
\begin{equation}\label{eqk}
\left\{
    \begin{array}{ll}
      k=\frac{\sqrt{Ey_1^2+4x_1}-2}{y_1}, & \hbox{if $ E \ge -\pi^2/4$;} \\
      k=-\frac{\sqrt{Ey_1^2+4x_1}+2}{y_1}, & \hbox{if $-\pi^2<E<-\pi^2/4$,}
    \end{array}
  \right.
\end{equation}
 and the optimal cost is
\begin{equation}\label{eC}
\left\{
  \begin{array}{ll}
    C=4\frac{\frac{E}{2}y_1-x_1-1+\sqrt{4x_1+g^{-1}\left(-\frac{y_1}{2\sqrt{x1}}\right)y_1^2}}{y_1}, & \hbox{if $ E \ge
-\pi^2/4$;} \\
    C=4\frac{\frac{E}{2}y_1-x_1-\sqrt{4x_1+g^{-1}\left(-\frac{y_1}{2\sqrt{x1}}\right)y_1^2}}{y_1}, & \hbox{if
$-\pi^2/4<E<-\pi^2/4$.}
  \end{array}
\right.
\end{equation}
We are now left to prove that, with the previous choice of $k$, one also has $x(0)=x_1$ and $x(t_1)=1$. With this goal,
it is sufficient to observe the following interesting geometric feature of solutions of \eqref{probl-2}: the quantity
$\lambda_1(s)x(s)+\lambda_2(s)y(s)$ is another constant of motion for \eqref{probl-2}, whose value set at $s=t_1$ is
$k$. Merging this information with \eqref{integrprimi}, we have
$$2c x(s)=E-(k-cy(s))^2$$
for all points $(x(s),y(s))$ of the solution of \eqref{probl-2}. In other terms, the trajectory $(x(s),y(s))$ always
belongs to the parabola
$$x(s)=-\frac{c}2 y^2(s)+ky(s)+1.$$
Then, when the trajectory reaches $y(0)=y_1$ and $t_1=2$, it holds
\begin{equation}
x(0)=\frac{k^2-E}{4}y_1^2+k y_1+1 = x_1,
\end{equation}
by plugging the explicit expression \eqref{eqk} of $k$.

Summing up, the optimal trajectory steering  $(x_1,y_1)$ to $(1,0)$ in time $t_1=2$ is the unique solution of
\eqref{probl-2} with final covector $(k, \frac{k^2-E}{2})$, where $k$ and $E$ are given by \eqref{eqk} and \eqref{e-E}.
We next prove Proposition \ref{prop-main} by applying the symmetry inverse transformations \eqref{tr-inv-psi} and
\eqref{dil-inv-psi}.

\medskip

\noindent{\sc Proof of Proposition \ref{prop-main}.} \
By \eqref{dil-inv-psi} with $r = \frac{t_1-t_0}{2}$ we find
$$
  \Psi(x_1,y_1,t_1;x_0,y_0,t_0) =
  \tfrac{2}{t_1-t_0} \Psi \left(\tfrac{x_1}{x_0},\tfrac{2(y_1-y_0)}{x_0(t_0-t_1)}, 2; 1,0,0\right).
$$
Moreover, the Hamiltonian of the optimal trajectory of \eqref{probl-2} corresponding to the right hand side of the
above equation is $\frac{(t_1-t_0)^2}{4} E$, where $E$ is the Hamiltonian of the optimal trajectory steering
$(x_1,y_1,t_1)$ to $(x_0,y_0,t_0)$. From \eqref{e-E} we obtain $\frac{(t_1-t_0)^2}{4} E =
g^{-1}\left(\frac{y_0-y_1}{(t_1-t_0)\sqrt{x_0 x_1}}\right)$, that gives \eqref{propE-prop}.  By using the first
expression in
\eqref{eC} of the $\Psi \left(\frac{x_1}{x_0},\frac{2(y_1-y_0)}{x_0(t_1-t_0)}, 2; 1,0,0\right)$, we obtain
\begin{equation*}
\begin{split}
  \Psi(x_1,y_1,t_1;x_0,y_0,t_0) = & \frac{2}{t_1-t_0} 4  \bigg( \frac{(t_1-t_0)^2}{4} \frac{E}{2}
\cdot \frac{2(y_1-y_0)}{x_0(t_0-t_1)} -\frac{x_1}{x_0}-1+\\
  & +\sqrt{4\frac{x_1}{x_0} + \frac{(t_1-t_0)^2}{4} E \left(\frac{2(y_1-y_0)}{x_0(t_1-t_0)}\right)^2 }\Bigg)
  \frac{x_0(t_1-t_0)}{2 (y_1-y_0)},
\end{split}
\end{equation*}
which, recalling that $y_0>y_1$, agrees with \eqref{V(E)-prop}. The proof of the second one is analogous.

\smallskip

In order to prove \eqref{asimpt1}, we claim that, for every $\eps \in ]0,1[$ there exists a positive $E_\eps$ such that for every $E > E_\eps$ it holds
\begin{equation}\label{bounds-E}
\frac{4}{(t_1-t_0)^2}\log^2\left(\frac{y_0-y_1}{(t_1-t_0)\sqrt{x_0x_1}} \right) < E <
\frac{4}{(1-\eps)^2(t_1-t_0)^2}\log^2\left(\frac{y_0-y_1}{(t_1-t_0)\sqrt{x_0x_1}} \right),
\end{equation}
where $E$ is the function defined in \eqref{e-E}. To prove the claim, we fix $\eps \in ]0,1[$
and we note that
\begin{equation}\label{sinh}
 \exp((1- \eps) x) < \frac{\sinh(x)}{x} < \exp(x),
\end{equation}
for every sufficiently large positive $x$. Recalling \eqref{propE-prop}, since $\tfrac{y_0-y_1}{(t_1-t_0)\sqrt{x_0x_1}}
\rightarrow +\infty$, we consider $g(r)$ in \eqref{func_g} with $r>0$. Then, from \eqref{sinh} it follows that
\begin{equation*}
  \exp\left(\frac{(1-\eps)(t_1-t_0) \sqrt{E}}{2}\right) < \frac{y_0-y_1}{(t_1-t_0)\sqrt{x_0x_1}}
  < \exp\left(\frac{(t_1-t_0) \sqrt{E}}{2}\right),
\end{equation*}
for any positive $E$ big enough. This proves \eqref{bounds-E}.
Moreover, for $E$ big enough, we have, for every arbitrary $\epsilon>0$
\begin{equation}\label{estimate}
0 \leq \frac{4x_1x_0}{(y_0-y_1)^2}=\frac{E}{\sinh^2 \big(\tfrac{(t_1-t_0)}{2}\sqrt{E}\big)}<\epsilon.
\end{equation}

We next consider the value function $\Psi$ as a function of $\frac{y_0-y_1}{(t_1-t_0)\sqrt{x_1x_0}}$. From the first
expression in
\eqref{V(E)-prop} and \eqref{estimate}, we obtain the following inequality
\begin{equation*} 
 \Psi(x_1,y_1,t_1;x_0,y_0,t_0) \le \frac{4}{(1-\eps)^2(t_1-t_0)} \log^2\left(\frac{y_0-y_1}{(t_0-t)\sqrt{x_0x_1}}
\right)+ \frac{4(x_1+x_0)}{y_0-y_1}
\end{equation*}
for every $E > E_\eps$. On the other hand, modifying if necessary the choice of $E_\eps$, we also have
\begin{equation*} 
  \Psi(x_1,y_1,t_1;x_0,y_0,t_0) \ge \frac{4(1-\eps)^2 }{(t_1-t_0)}
\log^2\left(\frac{y_0-y_1}{(t_1-t_0)\sqrt{x_0x_1}} \right)+ \frac{4(x_1+x_0)}{y_0-y_1}-2\epsilon
\end{equation*}
for every $E > E_\eps$.
This concludes the proof of \eqref{asimpt1}.

The proof of \eqref{asimpt2} is easier. It suffices to note that since, $\tfrac{y_0-y_1}{(t_1-t_0)\sqrt{x_1x_0}} \to
0$,
 we consider $g(r)$ in \eqref{func_g} with $r<0$, then $E\rightarrow
-\frac{4\pi^2}{(t_1-t_0)^2}$. From the second expression in \eqref{V(E)-prop} we have $$\lim_{E\rightarrow
-\frac{4\pi^2}{(t_1-t_0)^2}} \
\frac{\Psi(x_1,y_1,t_1;x_0,y_0,t_0)}{\frac{4(x_1+x_0)+4\sqrt{4x_1x_0}}{y_0-y_1}-\tfrac{4\pi^2}{(t_1-t_0)}}=1.$$
\hfill $\square$

\subsection{Lower bound in \eqref{e-twosidedbounds}}
In this section we give the proof of the lower bound in Theorem \ref{th-main} for a preliminary choice of the pole $z_0=(x_0,y_0,t_0)=(1,0,0)$. We pass to the general case at the end of Section \ref{sec6}.  We first prove the following

\begin{lemma} \label{lem-lwbd} There exists two positive constants $\kappa$ and $\varrho$, only depending on the
$L^\infty$ norms of $\widetilde a, \widetilde b$, and on $\inf \widetilde a$, such that
\begin{equation*} 
  \Gamma(1,- t, t; 1,0,0) \ge \frac{\kappa}{t^2}, \quad \text{for every} \quad  t \in ]0,  \varrho/4[.
\end{equation*}
\end{lemma}

\medskip
\noindent{\sc Proof.}\ We claim that, for every $r \in ]0, 1/2]$ we have
\begin{equation*} 
  \Gamma(x,y, t; \xi, \eta, \tau) \ge G_{r}(x,y, t;\xi, \eta, \tau),
\end{equation*}
for every $(x,y, t; \xi, \eta, \tau) \in \overline{\H0 (1,0,0)} \times \H0 (1,0,0)$, where $G_{r}(x,y, t;\xi, \eta,
\tau)$ is the Green function appearing in \eqref{eq-v-Green1}. The proof of Lemma \ref{lem-lwbd}
then follows from Lemma \ref{lem-bd-Green}.

In order to prove our claim, we fix $r \in ]0, 1/2]$. For every  non-negative $f \in C_0^\infty (\H0 (1,0,0))$ and for
every $(x,y, t) \in \overline{\H0 (1,0,0)}$ we set
\begin{align*} 
  & v_f (x,y,t) := \int_{\H0 (1,0,0)} G_r(x,y,t; \xi, \eta, \tau) f (\xi, \eta, \tau) d\x \, d \eta \, d \tau, \\
  & u_f (x,y,t) := \int_{\H0 (1,0,0)} \Gamma(x,y,t; \xi, \eta, \tau) f (\xi, \eta, \tau) d\x \, d \eta \, d \tau.
\end{align*}
Both $v_f$ and $u_f$ are solution of $\L u = -f$ in $\H0 (1,0,0)$. Moreover $u_f(x,y,t) \ge 0$ for every  $(x,y,t) \in
\partial (\H0 (1,0,0)) \cap \big\{t < r^2\big\}$. From \eqref{eq-v-Green2} and from the comparison principle we
then find $u_f \ge v_f$ in $\H0 (1,0,0)$. In other words, we have
\begin{equation*} 
  \int_{\H0 (1,0,0)} \big( \Gamma(x,y,t; \xi, \eta, \tau) - G_r(x,y,t; \xi, \eta, \tau) \big) f (\xi, \eta, \tau) d\x \,
d \eta \, d \tau \ge 0,
\end{equation*}
for every non-negative $f \in C_0^\infty (\H0 (1,0,0))$ and for every $(x,y, t) \in \overline{\H0 (1,0,0)}$. This proves
our claim.  \hfill $\square$

\medskip

We next state and prove the main result of this section.

\begin{proposition}\label{lower_bound_Fund}
Let $0<\eps<1$ be fixed arbitrarily. There exists a positive constant $c_{\eps,T}^-$ only depending on the operator
$\L$, on $\eps$ and on $T$ such that for every $(x,y,t) \in \R^+ \times \R \times]0,T]$ with $y< -\eps t$ it
holds
\begin{align}\label{lower_bound_gamma}
  &\G\left(x,y,t;1,0,0\right) \ge \frac{c_{\eps,T}^-}{t^2} 
   \exp \left(-C\Psi(x,y+\eps t, t-\eps t; 1,0,0) \right).
\end{align}
\end{proposition}

\medskip
\noindent{\sc Proof.}\ Let  $\eps \in ]0,1[$ be fixed, by Proposition \ref{p2a} and Lemma \ref{lem-lwbd} we have
\begin{align}
  \G(x,y,t;1,0,0) & \ge \eps^{-\beta} M^{-1-\frac{4(1-\eps) t}{\theta^2} -\frac{\Psi(x,y,t;1,-\eps t, \eps t)}{h}}
  \G(1,- \eps t, \eps t;1,0,0)\nonumber \\
 &\ge \eps^{-\beta}M^{-1-\frac{4T}{\theta^2}-\frac{\Psi(x,y,t;1,- \eps t, \eps t)}{h}} \frac{\kappa}{4 (\eps t)^2},
\label{first-step}
\end{align}
for every $(x,y,t) \in \R^+ \times \R \times ]0,T]$ with $y < - \eps t$. This proves  \eqref{lower_bound_gamma} for
$(x_0,y_0,t_0)=(1,0,0)$, with $c_{\eps,T}^-= \frac{\kappa}{ 4 \eps^{2+\beta}}M^{-1-\frac{4T}{\theta^2}}.$

\section{Upper Bound and Proof of the Main Theorem}\label{sec6}
\setcounter{section}{5} \setcounter{equation}{0} \setcounter{theorem}{0}

{In this section we prove the upper bound in
\eqref{e-twosidedbounds} for the fundamental solution of $\L$. For the scopes of this section  it is more convenient to write $\L$ in its divergence form \eqref{divergence_form}.

 To achieve the proof of Proposition \ref{th-ub}, we need to introduce some preliminary results on non-negative weak solutions $u$ to $ \L u = 0$ in $\R^+ \times \R \times ]T_0,T_1[$ and on non-negative weak solutions   $u$ to its formal adjoint $ \L^* u = 0$ in $\R^+ \times \R \times ]T_0,T_1[$.
 For this reason, we  consider operators with a \emph{zero order term}, namely
\begin{equation}\label{e-main-divform}
\L_1 u(x,y,t)=  - X^* \big( a X u \big) + (b-a) X u + c u + Y u.
\end{equation}
Clearly, $\L$ is the particular case of $\L_1$ that we obtain with $c=0$. With the the same notation, its formal adjoint
$\L_1^*$ is
 \begin{equation}\label{e-main-divform1}
\L_1^* u(x,y,t)=  - X^* \big( a X u \big) - X^* \big( (b-a) u \big) + c u - Y u.
\end{equation}
In the sequel we rely on the following assumption
\begin{equation} \label{eq-assumption-c}
  a, b, c, \p_x(xa), \p_x(xb) \quad \text{are  bounded and measurable functions.}
\end{equation}
Note that the same condition holds for $\L_1^*$. The existence of a fundamental solution for $\L$ is guaranteed if
we also suppose that the coefficients $a,b,c$ are smooth.

The main result of this section is the following
\begin{proposition}[Upper Bound]\label{th-ub}
Let $T_0,T_1$ be fixed and consider the  set $\R^+ \times \R \times ]T_0,T_1[$. Let $\L_1$ be the operator in
\eqref{e-main-divform}, and $\Gamma(x,y,t;1,0,0)$ be its fundamental solution. Denote by $M_1$ the
$L^\infty$-norm of $a(x,y,t)$ and $T=T_1-T_0$. Then, for every positive $\epsilon$, there exists a
positive constant $ C^+_{\epsilon}$, depending on the vector fields $X,Y$, on $\epsilon, T$ and on the
$L^\infty$-norm of $a(x,y,t)$ such that
\begin{equation}\label{upper-bond}
\Gamma(x,y,t; 1,0,0) \le \frac{C^+_{\epsilon}}{t^{2}}\exp
\left(- \tfrac{1}{16M_1}\Psi(x,y-\eps,t+\epsilon;1,0,0) \right)
\end{equation}
for every $(x,y,t) \in \R^+\times ]-\infty,0[ \times]0, T]$.
\end{proposition}
%
The proof of Theorem \ref{th-ub} is based on a local $L^\infty$ a
priori estimate for solution of $\L_1 u = 0$. In order to state precisely this estimate, we recall some notation. For every $(x_0,y_0,t_0) \in \R^+ \times \R^2$ and $r \in ]0,1[$ we consider the set $H_r(x_0,y_0,t_0)$ introduced in \eqref{spheres}.

\begin{proposition}\label{MoserL}
Let $(x_0,y_0,t_0)$ be any point of $\R^+ \times \R^2$, and let $r, \rho$ with $0<r/2 \le \rho <r < 1$. Let u be a
non-negative weak solution of $\L_1 u(x,y,t)=0$ in
$H_r(x_0,y_0,t_0) $ and let $u \in L^2(H_r(x_0,y_0,t_0))$. Then
\begin{equation}\label{Moser}
  \sup_{H_\rho(x_0,y_0,t_0)} u^p \le \frac{\bar{c}}{(r-\rho)^6}\int_{H_r(x_0,y_0,t_0)}u^p,
\end{equation}
where the constant $\bar{c}>0$ depends only on $\L_1, p$ and on the $L^\infty$ norm of $a, b, c$.
\end{proposition}

The proof of Proposition \ref{MoserL} relies on the analogous result proven in \cite[Theorem 1.4]{CPP} for the
Kolmogorov equation
with bounded coefficients. For the sake of simplicity we recall here its statement for a particular operator
strongly related to $\L_1$. For every $(x_0,y_0,t_0)$ and $r > 0$ we denote
$$
  \widetilde{H}_r(x_0,y_0,t_0):= \Big\{(x,y,t) \in \R^3  \ | \  |x-x_0|<r,\
  |y-y_0+x_0(t-t_0)|<r^3 , \ -r^2< t-t_0<0\Big\}.
$$
\emph{Let $\Omega$ be an open subset of $\R^3, (x,y,t) \in \Omega$ and consider $v(x,y,t)$ a
positive weak solution in $\Omega$ of the following equation}
\begin{equation} \label{kolmog1}
  \p_{x}(\widetilde{a}(x,y,t) \p_{x} v) + \widetilde {b}(x,y,t) \p_{x} v +x \p_{y} v+\widetilde {c}(x,y,t)v - \p_{t}
v=0.
\end{equation}
\emph{Assume that $\widetilde a$, $\widetilde b$ and $\widetilde {c}$ are measurable bounded continuous functions such
that $\inf_{\Omega} \widetilde a(x,y,t) >0$. Let $(x_0,y_0,t_0) \in \Omega$ and $\rho, r$ such that $0<r/2 \le \rho <r
\le 1$ and $\widetilde{H}_r(x_0,y_0,t_0) \subseteq \Omega$. Then, there exists a positive constant $c$
depending on the $L^\infty$ norm of $\widetilde a$, $\widetilde b$, $\widetilde {c}$ and on $p$ such that}
\begin{equation}\label{Moser-Kolmogorov}
  \sup_{\widetilde{H}_\rho(x_0,y_0,t_0)}  v^p \le \frac{c}{(r-\rho)^6}\int_{\widetilde{H}_r(x_0,y_0,t_0)} v^p.
\end{equation}
\emph{for every $u \in L^p(\widetilde{H}_r(x_0,y_0,t_0))$.}

\medskip

\noindent {\sc Proof of Proposition \ref{MoserL}.}\  We first note that $\L_1 u = 0$ reads as follows
\begin{equation}\label{rewriteL}
  \p_x(x^2a(x,y,t)\p_xu)+(b(x,y,t)-a(x,y,t))x\p_xu+c(x,y,t)u+x\p_yu-\p_tu=0
\end{equation}
so that it has the form \eqref{kolmog1}. Even if coefficents of $\L_1$ are unbounded and $\inf_{\R^+ \times \R^2} x^2 a
=
0$, estimate \eqref{Moser-Kolmogorov} holds on compact cylinders contained in $\R^+ \times \R^2$. However, we need to
show that
the constant $\bar c$ in \eqref{Moser} does not depend on $(x_0,y_0,t_0)$ and $r$.

We first fix $(x_0,y_0,t_0)=(1,0,0)$, so that the cylinders $H_r(1,0,0)$ and $\tilde{H}_r(1,0,0)$ coincide. We modify
the functions $a(x,y,t)$, $b(x,y,t)$ and $c(x,y,t)$ as we have done in Section \ref{sec3}
$$
  \widetilde{a}(x,y,t)=\varphi^2(x)a(x,y,t), \ \widetilde{b}(x,y,t)=
  \varphi(x)(b(x,y,t)-a(x,y,t)), \ \widetilde{c}(x,y,t)=\varphi(x)c(x,y,t)
$$
where $\varphi(x)$ is the function defined in $\eqref{e1a}$. Then the functions $\widetilde{a}, \widetilde{b}$ and
$\widetilde{c}$ are uniformly bounded, $\inf \widetilde{a}$ is strictly positive and \eqref{Moser-Kolmogorov} implies
\eqref{Moser} if $(x_0,y_0,t_0)=(1,0,0)$.

For a general $(x_0,y_0,t_0)$, we consider the function $w(x,y,t) := u \big((x_0,y_0,t_0) \circ (x,y,t) \big)$ and we
conclude the proof by the argument used in the proof of Proposition \ref{p1a}. \hfill $\square$

 \medskip

We next introduce a result that, combined with Proposition \ref{MoserL}, provides us with the asymptotic upper bound
of the fundamental solution of $\L_1$. We first introduce a suitable \emph{cut-off function}. Let choose $R>1$ and
consider the following function
\begin{equation}\label{cut-off}
  \chi_R(x,y)=g_R(x)h_R(y), \quad (x,y) \in \R^+ \times \R,
\end{equation}
where
\begin{description}
  \item[-] $ g_R(x)=\phi \left(\frac{\log^2(x)+1}{\log^2(R)+1} \right)$;
  \item[-] $\phi(s)$ is a continuous function such that $\phi(s)=1$ if $s \in [0,1/2]$ and $\phi(s)=0$ if $s \in
[1,+\infty[$;
  \item[-] $h(y)$ is a continuous function such that
  \begin{itemize}
    \item  $h(y)=1$ if $y \in [-R,R]$;
    \item  $h(y)=0$ if $ y \in ]-\infty, -R^2]\cup [R^2, +\infty[$;
    \item  $h(y)$ is a $C^2$ spline function with derivative bounded by $\tfrac{2}{R^2-R}$, if $ y \in
[-R^2;-R]\cup[R,R^2]$.
  \end{itemize}

\end{description}
We first observe that $ g_R(x)\neq 0$ only if  $ x \in [1/R,R]$ and
\begin{gather}\label{}
  |x \p_y \chi_R(x,y)| \le x |g_R(x)| |\p_yh_R(y)| \le \frac{2}{R-1}, \nonumber \\
  |x \p_x \chi_R(x,y)| \le x |h_R(y)| \|\phi^{\prime}\|_{L^{\infty}(\R)}\frac{2 \log(x)}{x(\log^2(R)+1)} \le
\|\phi^{\prime}\|_{L^{\infty}(\R)}\frac{2 \log(x)}{(\log^2(R)+1)}. \nonumber
\end{gather}
Therefore
\begin{gather}\label{}
|X \chi_R| \le C \frac{\log R}{\log^2 R +1} \rightarrow 0 \quad \mbox{as} \ R \rightarrow +\infty \nonumber \\
|Y \chi _R| \le |x \p_y \chi_R|\le \frac{2}{R-1} \rightarrow 0
\quad \mbox{as} \ R \rightarrow +\infty. \nonumber
\end{gather}

Now we are ready to state the following
\begin{proposition} \label{th-HJ}
Let $u \in L^2\left(\R^+ \times \R^2\right)$ be a weak solution of $\L_1 u=0$, and let $\Psi$ be the value function of
the control problem \eqref{cost-phi}. Then there exist two positive constants $m, M_1$ only depending on the $L^\infty$
norm of $a, \, b, \, c, \, x\p_xa, \, x\p_xb$, such that
\begin{equation}\label{costtheorem}
  \int_{\R^+\times \R}e^{-\frac{\Psi(x_1,y_1,s;x,y,t_1)}{8M_1}-mt_1}u^2(x,y, t_1)dx \, dy \le \int_{\R^+\times \R}
e^{-\frac{\Psi(x_1,y_1,s;x,y,t_0)}{8M_1}-mt_0}u^2(x,y,t_0)dx \, dy,
\end{equation}
for every $t_0, t_1$ with $t_0 < t_1$, and $(x_1,y_1,s) \in \R^+\times \R \times ]t_1,+\infty[$.
\end{proposition}

\medskip

\noindent{\sc Proof.}\ Fix $(x_1, y_1, t_1) \in \R^+ \times \R^2$, and $t_0 < t_1$, and recall that, for any
$(x_0,y_0,t_0) \in \R^3$, in view of \eqref{cost-phi} the function $(x,y,t) \mapsto \Psi(x_0,y_0,t_0; x,y,t)$ is a
classical solution of the Hamilton-Jacobi-Bellman equation (see \cite{Capuzzo})
$$
  Y\Psi+\frac{1}{4}(X \Psi)^2=0.
$$
We set $v(x,y,t):=\frac{1}{16M_1}\Psi(x_0,y_0,t_0; x,y,t)$ where $M_1$ is the
$L^\infty$-norm of $a$. Then $v$ satisfies
\begin{equation}\label{v_hamilton}
Yv+4M_1(Xv)^2=0.
\end{equation}
We prove \eqref{costtheorem} by showing that
\begin{equation}\label{proof_theorem}
  \lim_{R \rightarrow +\infty}\int_{{\R^+\times \R} \times [t_0,t_1]}\frac{d}{dt}\chi_R^2e^{-2v-mt}u^2 \le 0,
\end{equation}
 where $\chi_R$ is the cut-off function introduced above and the constant $m$ will be specified in the sequel. Let $u$
be a positive solution of $\L_1$ in the domain $\R^+ \times \R \times [t_0,t_1]$. We note that
$$\int_{{\R^+\times \R}  \times[t_0,t_1]}x\p_y\big(\chi_R^2e^{-2v-mt}u^2\big)=0$$
since the function $\chi_R(x,y)$ has compact support in $\R^+\times \R$. Therefore we obtain
\begin{align}
  &\int_{{\R^+\times \R} \times [t_0,t_1]}\frac{d}{dt}\chi_R^2e^{-2v-mt}u^2=-\int_{{\R^+\times \R}  \times
[t_0,t_1]}Y\big(\chi_R^2e^{-2v-mt}u^2\big) = \nonumber \\
 &=\int_{{\R^+\times \R} \times [t_0,t_1]}e^{-2v-mt}u^2 \left(-Y\big(\chi_R^2\big) + 2\chi_R^2Yv
-m\chi_R^2\right)-2\int_{{\R^+\times \R} \times[t_0,t_1]}\chi_R^2e^{-2v-mt}uYu.\label{last_term}
\end{align}
We first focus on the last term of \eqref{last_term}. By using the fact that $u$ is weak solution of $\L_1 u = 0$ one
gets
\begin{align}
  &A:=-2\int_{{\R^+\times \R} \times [t_0,t_1]}\chi_R^2e^{-2v-mt}uYu=-2\int_{{\R^+\times \R} \times
[t_0,t_1]}aX\big(\chi_R^2e^{-2v-mt}u\big)Xu+ \nonumber \\
  &2\int_{{\R^+\times \R} \times [t_0,t_1]}\big(\chi_R^2e^{-2v-mt}u\big)(b-a)Xu+2\int_{{\R^+\times \R} \times
[t_0,t_1]}c\chi_R^2e^{-2v-mt}u^2=: A_1+A_2+A_3. \label{first_term}
\end{align}
Consider the first term in \eqref{first_term} and compute the derivatives
\begin{align}
  &A_1=-2\int_{{\R^+\times \R} \times [t_0,t_1]}aX\big(\chi_R^2e^{-2v-mt}u\big)Xu =-4\int_{{\R^+\times \R} \times
[t_0,t_1]}a\chi_Re^{-2v-mt}uXuX\chi_R+ \nonumber\\
 &4\int_{{\R^+\times \R} \times [t_0,t_1]}a\chi_R^2e^{-2v-mt}uXuXv-2\int_{{\R^+\times \R} \times
[t_0,t_1]}a\chi_R^2e^{-2v-mt}(Xu)^2=: B_1+B_2+B_3. \label{formula_1}
\end{align}
By using Young inequality, it follows
\begin{align}
  &B_1=-4\int_{{\R^+\times \R} \times
[t_0,t_1]}a\chi_Re^{-2v-mt}uXuX\chi_R \leq 4\int_{{\R^+\times \R} \times [t_0,t_1]}a\chi_Re^{-2v-mt}\, |Xu| \,
|uX\chi_R| \le \nonumber  \\
  &\int_{{\R^+\times \R} \times [t_0,t_1]}a\chi_R^2e^{-2v-mt}(Xu)^2+4\int_{{\R^+\times \R} \times
[t_0,t_1]}ae^{-2v-mt}u^2(X\chi_R)^2=:
  C_1+C_2, \label{young}
\end{align}
Merging the inequalities  \eqref{formula_1} and \eqref{young}, since $B_3=-2C_1$, we conclude
\begin{align}
  &A_1=-\int_{{\R^+\times \R} \times [t_0,t_1]}a\chi_R^2e^{-2v-mt}(Xu)^2+
  B_2+C_2\leq  \nonumber\\
  &4\int_{{\R^+\times \R} \times [t_0,t_1]}a\chi_R^2e^{-2v-mt}u^2(Xv)^2+C_2\leq 4M_1\int_{{\R^+\times \R} \times
[t_0,t_1]}\chi_R^2e^{-2v-mt}u^2(Xv)^2+C_2. \label{first_inequality}
\end{align}
Now consider the second term in \eqref{first_term}. Start from integration by parts formula
$$A_2=2\int_{{\R^+\times \R} \times [t_0,t_1]}u X^*\left((b-a)\big(\chi_R^2e^{-2v-mt}u\big)\right).$$
Reminding that $X^*=-X-1$, similarly to \eqref{formula_1}, \eqref{young} and \eqref{first_inequality} we have
\begin{align}
&A_2 \le -\int_{{\R^+\times \R} \times[t_0,t_1]}(b-a)\chi_R^2e^{-2v-mt}u^2+
 4M_1\int_{{\R^+\times \R} \times [t_0,t_1]}\chi_R^2e^{-2v-mt}u^2(Xv)^2+ \nonumber\\
 &\frac{1}{4M_1}\int_{{\R^+\times \R} \times [t_0,t_1]}(b-a)^2\chi_R^2e^{-2v-mt}u^2+
  \int_{{\R^+\times \R} \times [t_0,t_1]}X\big(\chi_R^2\big)e^{-2v-mt}u^2+  \nonumber\\
  &\int_{{\R^+\times \R} \times [t_0,t_1]}X(b-a)\chi_R^2e^{-2v-mt}u^2 \le 4M_1\int_{{\R^+\times \R} \times
[t_0,t_1]}\chi_R^2e^{-2v-mt}u^2(Xv)^2+ \nonumber \\
  &\int_{{\R^+\times \R} \times [t_0,t_1]}X\big(\chi_R^2\big)e^{-2v-mt}u^2+\frac{3}{4}m\int_{{\R^+\times \R} \times
[t_0,t_1]}\chi_R^2e^{-2v-mt}u^2. \label{second_inequality}
\end{align}
by setting
\begin{align}
  m := 4\, \max \Big\{\tfrac{1}{4M_1} \|(b-a)^2 & \|_{L^\infty({\R^+\times \R}\times[t_0,t_1])},\,
  \|X(b-a)\|_{L^\infty({\R^+\times \R}\times[t_0,t_1])}, \, \nonumber  \\
   &\quad 2\| c \|_{L^\infty({\R^+\times \R}\times[t_0,t_1])},
  \, \|b-a\|_{L^\infty({\R^+\times \R}\times[t_0,t_1])}\Big\}. \label{defn_m}
\end{align}
Going back to \eqref{first_term}, combining \eqref{first_inequality}, \eqref{second_inequality} and $$A_3 \le
\frac{1}{4}m\int_{{\R^+\times \R} \times [t_0,t_1]}\chi_R^2e^{-2v-mt}u^2,$$ we have
\begin{align}
&A \le 8M_1\int_{{\R^+\times \R} \times
[t_0,t_1]}\chi_R^2e^{-2v-mt}u^2(Xv)^2+
   4\int_{{\R^+\times \R} \times [t_0,t_1]}ae^{-2v-mt}u^2(X\chi_R)^2+  \nonumber \\
   & \int_{{\R^+\times \R} \times [t_0,t_1]} X(\chi_R^2)e^{-2v-mt}u^2+
 m \int_{{\R^+\times \R} \times [t_0,t_1]}\chi_R^2e^{-2v-mt}u^2. \label{main_ineq}
\end{align}
Merging \eqref{main_ineq} with \eqref{first_term} ad \eqref{last_term} we conclude that
\begin{align*}
&\int_{{\R^+\times \R} \times [t_0,t_1]}\frac{d}{dt}\chi_R^2e^{-2v-mt}u^2\leq 2\int_{{\R^+\times \R} \times
[t_0,t_1]}\big(\chi_R^2e^{-2v-mt}u^2 \big) \big[Yv+4M_1(Xv)^2 \big]+\\
 &\int_{{\R^+\times \R} \times [t_0,t_1]}\big(e^{-2v-mt}u^2 \big)\big(-Y\big(\chi_R^2\big)+4a(X\chi_R)^2+X\big(\chi_R^2
\big)\big).
\end{align*}
The first integral is zero since $v$ satisfies the Hamilton-Jacobi-Bellman equation \eqref{v_hamilton}, and
\eqref{costtheorem} simply follows by letting $R\rightarrow +\infty$. \hfill $\square$

\medskip

The next Lemma is crucial to prove Theorem \ref{th-ub}.

\begin{lemma}\label{lemma1}
Let $\epsilon$ be a fixed positive constants. Then there exists a constant $c_{\epsilon}>0$, only
depending on $\L_1$ and  $\epsilon$ such that
\begin{equation}\label{inequality1}
u^2\big(1, 0,t/2\big) \le c_{\epsilon} \int_{\R^+\times
\R}e^{-\frac{1}{8M}\Psi\big(1,-\eps,\frac{t}{2} +\epsilon, \xi,\eta,0\big)}u^2(\xi,\eta,0)d\xi d\eta
\end{equation}
for every non-negative weak solution $u$ of $\L_1 u=0$ in $\R^+ \times \R \times ]T_0,T_1[$.
\end{lemma}
\medskip

\noindent {\sc Proof.} Let $\eps>0$ be fixed and let $r \in ]0, 1[$ be such that $r^3< \eps$. By
Proposition \ref{MoserL}, with $p=2$,  we have
\begin{equation}\label{identity_1Moser}
  u^2\big(1,0,t/2\big) \le \sup_{H_{r/2}^-\big(1,0,t/2\big)}u^2(\xi, \eta, \tau) \le
\frac{c}{(r/2)^6}\int_{H_{r}\big(1,0,t/2\big)}u^2(\xi, \eta, \tau)d\xi d \eta d \tau
\end{equation}
for every $t \in ]T_0,T_1[$.  Multiply and divide the integrand of the above inequality by the quantity
\begin{equation}\label{function_estimate}
e^{\tfrac{1}{8M_1}\Psi\big(1,-\eps,\frac{t}{2} +\epsilon;\xi,\eta,\tau\big)+m\tau}.
\end{equation} Note that, as $r^3 < \eps$, the function $(\xi, \eta, \tau) \mapsto \Psi\big(1,-\eps,\frac{t}{2}
+\epsilon;\xi,\eta,\tau\big)$ is well defined, continuous and bounded in the set $\overline{H_{r/2}\big(1,0,t/2\big)}$.

Therefore, we denote by $C_{\epsilon}$ the maximum of the function in \eqref{function_estimate} in the set
$\overline{H_{r/2}\big(1,0,t/2\big)}$, which is uniform with respect to $t \in ]T_0,T_1[$. We then find

  \begin{align*}
   u^2\big(1,0,t/2\big) &\le \frac{c}{(r/2)^6} \cdot C_{\epsilon} \int_{H_{r}\big(1,0,t/2\big)}
   e^{-\tfrac{1}{8M_1}\Psi\big(1,-\eps,\frac{t}{2} + \epsilon;\xi,\eta,\tau\big)-m\tau}u^2(\xi,\eta,\tau)d\xi d \eta d
\tau \nonumber \\
  &\le \frac{c}{(r/2)^6} \cdot C_{\epsilon}\int_{t/2-r^2}^{t/2}\int_{\R^+\times
\R}e^{-\tfrac{1}{8M_1}\Psi\big(1,-\eps,\frac{t}{2} +\epsilon;\xi,\eta,\tau\big)-m\tau}u^2(\xi,\eta,\tau)d\xi d \eta d
\tau \nonumber \\
& \mbox{(by Proposition \ref{th-HJ}, with $t_0=0$ and $t_1 = \tau$)}\\
  &\le \frac{c}{(r/2)^6} \cdot C_{\epsilon}\int_{t/2-r^2}^{t/2}\int_{\R^+\times
\R}e^{-\tfrac{1}{8M_1}\Psi\big(1,-\eps,\frac{t}{2} +\epsilon;\xi,\eta,0\big)}u^2(\xi,\eta,0)d\xi d \eta \nonumber \\
  & \le \frac{c}{(r/2)^6} \cdot C_{\epsilon}\int_{\R^+\times\R}
  e^{-\tfrac{1}{8M_1}\Psi\big(1,-\eps,\frac{t}{2} + \epsilon;\xi,\eta,0\big)}u^2(\xi,\eta,0)d\xi d \eta.
\end{align*}
which gives \eqref{inequality1} by setting $c_{\epsilon} := \frac{c}{(r/2)^6} \cdot C_{\epsilon}$.
 \hfill $\square$

 \medskip

We finally introduce a last result we need to prove Proposition \ref{th-ub}.
The following Proposition is a direct consequence of Proposition \ref{MoserL} and involves the fundamental solution $\Gamma$ of $\L$.

\begin{proposition}\label{results}
Let $\Gamma$ be a fundamental solution of $\L_1$ and fix $(x,y,t), \ (x_0,y_0,t_0)$ in $\R^+\times \R^2$ with
$y<y_0$ and $T_0\le t_0<t \le T_1$. Define $T=T_1-T_0$. Then, there exist a positive constant $C_T$ depending on
the operator $\L_1$ and on $T$ such that  the following upper bounds hold for $\Gamma$
\begin{description}
  \item[i)] $\Gamma(x,y,t;x_0,y_0,t_0)\le \frac{C_T}{(t-t_0)^{2}}$ ;
  \item[ii)] $\int_{\R^+\times \R }\Gamma^2(x,y,t;x_0,y_0,t_0)dx_0dy_0 \le \frac{C_T}{(t-t_0)^{2}} ;$
\end{description}
\end{proposition}

\medskip

\noindent {\sc Proof.} We only prove \textbf{i)}, since \textbf{ii)} is its direct consequence reminding that
$$\int_{\R^+\times \R }\Gamma(x,y,t;x_0,y_0,t_0)dx_0dy_0=1.$$
We first fix $0<t-t_0<1$ and, by using Proposition \ref{MoserL}, we have
\begin{align}
\Gamma(x,y,t;x_0,y_0,t_0) &\le \sup_{H_{\sqrt{t-t_0}/2}(x,y,t)} \Gamma(\cdot, \cdot, \cdot;x_0,y_0,t_0) \nonumber \\
 &\le \frac{\bar{C}}{(t-t_0)^3} \int_{{H_{\sqrt{t-t_0}}}(x,y,t)}\Gamma(\xi, \eta, \tau;x_0,y_0,t_0)d \xi d \eta d \tau
\nonumber \\
 & \le \frac{\bar{C}}{(t-t_0)^3} \int_{t-(t-t_0)}^t d \tau \int_{\R^+ \times \R}\Gamma(\xi, \eta, \tau;x_0,y_0,t_0)d
\xi d \eta = \frac{\bar{C}}{(t-t_0)^2} \label{firstineq}
\end{align}
since $\int_{\R^+ \times \R}\Gamma(\xi, \eta, \tau;x_0,y_0,t_0)d \xi d \eta $ is finite in view of \eqref{integrals}. If $t-t_0 \ge 1$ we set $
\nu=\frac{t-t_0}{T}<1$, and starting from the reproduction property we have
\begin{align*}
\Gamma(x,y,t;x_0,y_0,t_0)&=\int_{\R^+ \times \R}\Gamma(x,y,t;\xi, \eta, t_0+\nu)\Gamma(\xi, \eta,t_0+\nu;x_0,y_0,t_0)d
\xi d\eta \\
&\le \frac{C_T}{(t-t_0)^2}\int_{\R^+ \times \R}\Gamma(x,y,t;\xi, \eta, t_0+\nu)d \xi d \eta \le\frac{C_T}{(t-t_0)^2}
\end{align*}
by \eqref{firstineq} where $C_T=\bar{C}\,T^2$ and $\int_{\R^+ \times \R}\Gamma(x,y,t;\xi, \eta, t_0+\nu) d \xi d
\eta=1$.
  \hfill $\square$

\medskip

We are now ready to prove the main proposition of this Section.

\medskip

\noindent {\sc Proof of Proposition \ref{th-ub}.} \
 Let $\eps >0$ be fixed and let $\Gamma(x,y,t;1,0,0)$ be the fundamental solution of $\L_1$ and
$(x,y,t) \in \R^+ \times ]-\infty, 0[\times]T_0,T_1[$. We define
\begin{align*}
  D_1 & = \, \big\{(\xi, \eta) \in \R^+ \times \R^- \, |\  \Psi\big(x,y-\eps,t+\epsilon/2, \xi, \eta, t/2\big) \leq
\Psi(\xi, \eta, t/2; 1,0,-\epsilon/2) \big\}, \\
   D_2 & = \, \big\{(\xi, \eta) \in \R^+ \times \R^- \, | \Psi(x,y-\eps,t+\epsilon/2, \xi, \eta, t/2) > \Psi(\xi,
\eta, t/2; 1,0,-\epsilon/2) \big\},
\end{align*}
Starting from the reproduction property of $\Gamma$
\begin{align*}
\Gamma(x,y,t;1,0,0) & = \int_{\R^+\times \R^-}\Gamma(x,t,y; \xi, \eta, t/2) \Gamma(\xi, \eta,t/2, 1,0,0) d \xi d \eta \\
 & = \int_{D_1}\Gamma(x,t,y; \xi, \eta, t/2) \Gamma(\xi, \eta,t/2, 1,0,0) d \xi d \eta+ \\
 & +\int_{D_2}\Gamma(x,t,y; \xi, \eta, t/2) \Gamma(\xi, \eta,t/2, 1,0,0) d \xi d \eta \\
 & \leq \frac{C_T}{t^{2}} \left(\left(\int_{D_1} \Gamma^2(\xi, \eta,t/2, 1,0,0) d \xi d \eta\right)^{\tfrac{1}{2}}
+\left(\int_{D_2} \Gamma^2(x,y,t;\xi, \eta,t/2) d \xi d \eta\right) ^{\tfrac{1}{2}} \right)
\end{align*}
 where $T=T_1-T_0$ and the last inequality follows from H\"{o}lder inequality and \eqref{firstineq}.\\
 We now introduce the sets
\begin{align*}
 \tilde{D}_1 & = \, \big\{(\xi, \eta) \in \R^+ \times \R^- \, |\
\Psi\big(x,y-\eps,t+\epsilon/2,1,0,-\epsilon/2\big) \leq 2\Psi(\xi, \eta,t/2, 1,0,-\epsilon/2)\big\}, \\
   \tilde{D}_2 & = \, \big\{(\xi, \eta) \in \R^+ \times \R^- \, |\  \Psi(x,y-\eps,t+\epsilon/2, 1,0,-\epsilon/2)
\leq 2\Psi(x,y-\eps, t+\epsilon/2;\xi, \eta,t/2) \big\},
\end{align*}
and we note that  $D_1 \subseteq \tilde{D}_1$ and $D_2 \subseteq \tilde{D}_2$ as a consequence of the \emph{triangular
inequality} of the value function:
$$\Psi(x_0,y_0,t_0;x,y,t) \le \Psi(x_0,y_0,t_0;\xi,\eta,\tau)+\Psi(\xi,\eta,\tau;x,y,t)$$
for arbitrary points $(x_0,y_0,t_0), (\xi,\eta,\tau), (x,y,t)$ belonging to $\R^+\times \R \times ]T_0,T_1[$ with
$y > \eta> y_0$ and $T_0\le t<\tau < t_0 \le T_1$. Hence
$$
\Gamma(x,y,t;1,0,0) \leq \frac{C_T}{t^{2}} \left( \left(\int_{\tilde{D}_1} \Gamma^2(\xi, \eta,t/2, 1,0,0) d \xi d \eta
\right)^{\tfrac{1}{2}} + \left(\int_{\tilde{D}_2} \Gamma^2(x,y,t;\xi, \eta,t/2) d \xi d \eta
\right)^{\tfrac{1}{2}}\right)
$$
We now claim that:

 \begin{align}
  &\int_{\tilde{D}_1} \Gamma^2(\xi, \eta,t/2, 1,0,0) d \xi d \eta \le c_{\epsilon}
e^{-\tfrac{1}{16M_1}\Psi(x,y-\eps,t+\epsilon/2; 1,0,-\epsilon/2)} \label{upper_bound_1}\\
  &\int_{\tilde{D}_2} \Gamma^2(x,y,t;\xi, \eta,t/2) d \xi d \eta  \le c_{\epsilon}
e^{-\tfrac{1}{16M_1}\Psi(x,y-\eps,t+\epsilon/2; 1,0,-\epsilon/2)} \label{upper_bound_2}
 \end{align}
 where $c_{\epsilon}$ is a positive constant depending on $\L$ and $\epsilon$. We first prove
\eqref{upper_bound_2} and we define the functions
 $$ v(z,w,s)=\int_{\tilde{D}_2}\Gamma(z,w,s;\xi, \eta,t/2)\Gamma(x,y,t;\xi, \eta,t/2) d \xi d \eta, \quad
u(z,w,s)=v((x,y,t/2) \circ (z,w,s))
 $$
 We further note that the functions $u$ and $v$ satisfy the following properties:
 \begin{description}
   \item[i)] $v(z,w,s)$ is a solution of $\L_1 v(z,w,s)=0$ in $\R^+\times \R \times [t/2, T_1[$. Then $u(z,w,s)$ is a
solution of
$\L_{1_{\bar{z}}}u(z,w,s)= 0$ in $\R^+ \times \R \times ]0,T_1[$ where $\bar{z}=(x,y,t/2)$ and
\begin{eqnarray}
 \L_{1_{\bar{z}}}u(z,w,s)&=& z \partial_{z} \big( a(xz,y+xw,t/2+s) z \partial_{z} u \big) + z \, b
(xz,y+xw,t/2+s)\partial_z u+\nonumber\\
& & + z \partial_{w}u+ c(xz,y+xw,t/2+s)u - \partial_{t} u.\label{L_modified}
 \end{eqnarray}

   \item[ii)] the function $v$ satisfies the initial condition $v(z,w,t/2)=\Gamma(x,y,t;z,w,t/2)
\mathbf{1}_{\tilde{D}_2}(z,w)$;
   \item[iii)]it holds $u(1,0,t/2)=\int_{\tilde{D}_2}\Gamma^2(x,y,t;\xi,\eta,t/2)d \xi d \eta$.
\end{description}
where $\mathbf{1}_{\tilde{D}_2}(z,w)$ denotes the characteristic function of the set $\tilde{D}_2$. In virtue of Lemma
\ref{lemma1} we have
$$ u^2\big(1,0,t/2\big) \le c_{\epsilon} \int_{\R^+\times
\R}e^{-\frac{1}{8M_1}\Psi\big(1,-\eps,\tfrac{t}{2}+\epsilon/2, z,w,0\big)}u^2(z,w,0)dzdw$$
By observing that $\Psi\big(1,-\eps,\tfrac{t}{2}+\epsilon/2, z,w,0\big)=\Psi\big(x,y-\eps,t+\epsilon/2;x,y,t/2)
\circ (z,w,0)\big)$, by the change of variable $ (\xi,\eta,t/2) = (x,y,t/2) \circ (z,w,0)$ and by properties
\textbf{ii)}
and \textbf{iii)}, we get
\begin{align*}
 \left(\int_{\tilde{D}_2}\Gamma^2(x,y,t;\xi, \eta, t/2) d \xi d \eta \right)^2 &= u^2\big(1,0,t/2\big)\\
  &\le c_{\epsilon} \int_{\tilde{D}_2}e^{-\frac{1}{8M_1}\Psi\big(x,y-\eps,t+\epsilon/2;
\xi,\eta,t/2\big)}\Gamma^2(x,y,t;\xi,\eta,t/2)d \xi d \eta
\end{align*}
We finally obtain \eqref{upper_bound_2} by recalling the definition of $\tilde{D}_2$
$$\left(\int_{\tilde{D}_2}\Gamma^2(x,y,t;\xi, \eta, t/2) d \xi d \eta \right)^2 \le c_{\epsilon}
e^{-\frac{1}{16M_1}\Psi\big(x,y-\eps,t+\epsilon/2,1,0,-\epsilon/2\big)}\int_{\tilde{D}_2}\Gamma^2(x,y,t;\xi,\eta,
t/2)d \xi d \eta$$
and the result immediately follows by dividing by $\int_{\tilde{D}_2}\Gamma^2(x,y,t;\xi, \eta, t/2) d \xi d \eta$ and by
recalling that $\Psi\big(x,y-\eps,t+\epsilon/2,1,0,-\epsilon/2\big)=\Psi\big(x,y-\eps,t+\epsilon,1,0,0)$

\smallskip

The proof of inequality \eqref{upper_bound_1} is analogous to \eqref{upper_bound_2}. Indeed,
consider the function
$$ v_2(z,w,s)=\int_{\tilde{D}_1}\Gamma(\xi, \eta,t/2;z,w,s)\Gamma(\xi, \eta,t/2;1,0,0) d \xi d \eta,
 $$
which is a non-negative solution to $\L_1^*v_2 = 0$  with final data $v_2(z,w,t/2) = \Gamma(z,w,t/2;1, 0, 0)$  if $(z,w)
\in \tilde{D}_1$ and $ v_2(z,w,t/2) =0$ if $(z,w) \notin \tilde{D}_1$. Notice that the coefficients of $\L_1^*$
satisfy the same assumptions \eqref{e-assumption1} and \eqref{eq-assumption-c} made on $\L_1$, then all the properties
shown for the function $(x,y,t) \mapsto \G(x,y,t;\xi,\eta,\tau)$ and used to prove \eqref{upper_bound_2}, also hold
for $(x,y,t) \mapsto \G(\xi,\eta,\tau;x,y,t)$ (which is the  fundamental solution of $\L_1^*u=0$) and they can be used
to prove \eqref{upper_bound_1}. This proves the claim. \hfill $\square$

\medskip

We are now ready to prove the main result of our article.

\medskip

\noindent {\sc Proof of Theorem \ref{th-main}}. Let $\G(x,y,t;x_0,y_0,t_0)$ denote the fundamental solution of $\L$ in
\eqref{e-main} and $(x,y,t), (x_0,y_0,t_0)$ in $\R^+\times \R\times [0,T]$ with $y<y_0$ and $t>t_0$. If
$(x_0,y_0,t_0)=(1,0,0)$, the lower bound of $\G$ follows from Proposition \ref{lower_bound_Fund},  whereas the upper
bound follows from Proposition \ref{th-ub}.
For a general choice of $z_0 = (x_0, y_0, t_0)$ it suffices to note that the function
\begin{equation}\label{trasl_sol_fund}
  \Gamma_{z_0}(x,y,t;1,0,0)= x_0^2 \, \Gamma((x_0,y_0,t_0) \circ (x,y,t);x_0,y_0,t_0)
\end{equation}
is the fundamental solution of the operator $\L_{z_0}$ defined in \eqref{Kaz0}, with singularity at $(1,0,0)$. As
noticed in Remark \ref{rem-invariance}, it satisfies assumptions \eqref{e-assumption1} with
the same constants $\lambda$  used for $\L$, then \eqref{lower_bound_gamma} and \eqref{upper-bond}
applies to $\Gamma_{z_0}$. If one consider the lower estimates, we find
\begin{equation*}
  \Gamma((x_0,y_0,t_0) \circ (x,y,t);x_0,y_0,t_0) \ge  \frac{c_{\eps,T}^-}{x_0^2 t^2}
\exp \left(-C^-\Psi(x,y,t; 1,- \eps t, \eps t)\right),
\end{equation*}
that can be written equivalently as follows
\begin{equation*}
  \Gamma((x,y,t;x_0,y_0,t_0) \ge  \frac{c_{\eps,T}^-}{x_0^2 (t-t_0)^2}
\exp \left(-C^-\Psi((x_0,y_0,t_0)^{-1} \circ (x,y,t); 1,- \eps (t-t_0), \eps (t-t_0))\right).
\end{equation*}
The conclusion follows by applying the invariance property \eqref{transl-phi} of $\Psi$:
\begin{align*}
\Psi((x_0,y_0,t_0)^{-1} \! \circ \! (x,y,t);1,-\eps (t-t_0), \eps(t-t_0))   \!
& = \Psi(x,y,t;(x_0,y_0,t_0)  \! \circ  \! (1, -\eps (t-t_0), \eps(t-t_0))) \\
&=\Psi(x,y,t;x_0,y_0-\eps (t-t_0) x_0,t_0+\eps(t-t_0))\\
&=\Psi(x,y+\eps (t-t_0) x_0,t-\eps(t-t_0);x_0,y_0,t_0).
\end{align*}

The proof of the upper bound is analogous. \hfill $\square$


\def\cprime{$'$} \def\cprime{$'$} \def\cprime{$'$}
  \def\lfhook#1{\setbox0=\hbox{#1}{\ooalign{\hidewidth
  \lower1.5ex\hbox{'}\hidewidth\crcr\unhbox0}}} \def\cprime{$'$}
  \def\cprime{$'$} \def\cprime{$'$}

\end{document}